\documentclass[11pt]{amsart}
\usepackage{amssymb,mathrsfs}
\newtheorem{theorem}{Theorem}[section] 
\newtheorem{lemma}[theorem]{Lemma} 
\newtheorem{proposition}[theorem]{Proposition} 
\newtheorem{corollary}[theorem]{Corollary} 
\theoremstyle{definition}
\newtheorem{definition}[theorem]{Definition}
\newtheorem{remark}[theorem]{Remark}

\begin{document}

\title{Systolic inequalities and the Horowitz-Myers conjecture}
\author{Simon Brendle and Pei-Ken Hung}
\address{Columbia University \\ 2990 Broadway \\ New York NY 10027 \\ USA}
\address{University of Illinois \\ 1409 W. Green Street \\ Urbana IL 61801 \\ USA}
\thanks{The authors are grateful to Professors Aghil Alaee, Piotr Chru\'sciel, Nick Edelen, Gerhard Huisken, Christos Mantoulidis, and Leon Simon for discussions and helpful comments. The first author was supported by the National Science Foundation under grant DMS-2403981 and by the Simons Foundation. He acknowledges the hospitality of T\"ubingen University and the Centre de Recherches Math\'ematiques, Montr\'eal, where part of this work was carried out.}
\begin{abstract}
Let $n$ be an integer with $3 \leq n \leq 7$, let $M$ be a compact manifold of dimension $n$ with boundary $\partial M$, and let $g$ be a Riemannian metric on $M$ with scalar curvature at least $-n(n-1)$. Under a topological assumption on $M$, we establish an inequality relating the infimum of the boundary mean curvature to the systole of the boundary $\partial M$. As a consequence, we obtain a new positive energy theorem, with equality being attained by the Horowitz-Myers metrics. 
\end{abstract}
\maketitle 

\section{Introduction} 

In this paper, we establish the following geometric inequality for two-dimensional surfaces.

\begin{theorem} 
\label{2D.inequality}
Let us fix a real number $N > 2$. Let $\Sigma$ be a compact, connected, orientable surface with non-empty boundary $\partial \Sigma$, and let $g$ be a Riemannian metric on $\Sigma$. We denote by $K$ the Gaussian curvature of $\Sigma$, by $\kappa$ the geodesic curvature of the boundary $\partial \Sigma$, and by $\eta$ the outward-pointing unit normal vector field to $\partial \Sigma$. Let $\psi$ be a smooth function on $\Sigma$ such that 
\[-2 \, \Delta \psi - \frac{N-1}{N-2} \, |\nabla \psi|^2 + 2K + N(N-1) \geq 0\] 
at each point in $\Sigma$. If $\partial \Sigma$ is connected, then 
\[2 \, |\partial \Sigma|^N \, \inf_{\partial \Sigma} (\langle \nabla \psi,\eta \rangle + \kappa - (N-1)) < \Big ( \frac{4\pi}{N} \Big )^N.\] 
If $\partial \Sigma$ is disconnected, then 
\[\inf_{\partial \Sigma} (\langle \nabla \psi,\eta \rangle + \kappa - (N-1)) \leq 0.\] 
\end{theorem}

To prove Theorem \ref{2D.inequality}, we break the discussion into two cases. If $\partial \Sigma$ is connected, the proof relies on a monotonicity formula (see Corollary \ref{monotonicity} below). This argument shares some common features with the proof of the fill radius estimate in the groundbreaking work of Gromov and Lawson (see \cite{Gromov-Lawson}, Section 10). If $\partial \Sigma$ is disconnected, we consider curves joining two different connected components of $\partial \Sigma$, and minimize a weighted length functional. The latter argument is inspired by the work of Schoen and Yau \cite{Schoen-Yau2}. 

In higher dimensions, we prove the following estimate.

\begin{theorem} 
\label{higher.dim.inequality.1}
Let us fix an integer $n$ with $3 \leq n \leq 7$ and a real number $N > n$. Let $M$ be a compact, connected, orientable manifold of dimension $n$ with non-empty boundary $\partial M$. Suppose that $\theta_0$ is a smooth map from $\partial M$ to $S^1$, and $(\theta_1,\hdots,\theta_{n-2})$ is a smooth map from $M$ to $T^{n-2}$. We assume that the map $(\theta_0,\theta_1|_{\partial M},\hdots,\theta_{n-2}|_{\partial M}): \partial M \to S^1 \times T^{n-2}$ has non-zero degree. We denote by $\Theta_0$ the pull-back of the volume form on $S^1$ under the map $\theta_0: \partial M \to S^1$. Note that $\Theta_0$ is a closed one-form on $\partial M$. Let $g$ be a Riemannian metric on $M$. We denote by $R_M$ the scalar curvature of $M$, by $H_{\partial M}$ the mean curvature of the boundary $\partial M$, and by $\eta$ the outward-pointing unit normal vector field to $\partial M$. Let $\varphi$ be a smooth function on $M$ such that 
\[-2 \, \Delta_M \varphi - \frac{N-n+1}{N-n} \, |\nabla^M \varphi|^2 + R_M + N(N-1) \geq 0\] 
at each point in $M$. Then 
\[2 \, \sigma^N \, \inf_{\partial M} (\langle \nabla^M \varphi,\eta \rangle + H_{\partial M} - (N-1)) < \Big ( \frac{4\pi}{N} \Big )^N,\] 
where $\sigma$ denotes the length of the shortest closed curve $\alpha$ in $\partial M$ satisfying $\int_\alpha \Theta_0 \neq 0$.
\end{theorem}

In particular, if $M$ is diffeomorphic to $B^2 \times T^{n-2}$, then the topological assumptions in Theorem \ref{higher.dim.inequality.1} are satisfied. To deduce Theorem \ref{higher.dim.inequality.1} from Theorem \ref{2D.inequality}, we construct a minimal slicing with free boundary. The minimal slicing technique was introduced in the fundamental work of Schoen and Yau \cite{Schoen-Yau1}. This argument is closely related to the torical symmetrization procedure developed by Gromov and Lawson \cite{Gromov-Lawson} (see also \cite{Fischer-Colbrie-Schoen} and \cite{Gromov}).

If we put $\varphi=0$ in Theorem \ref{higher.dim.inequality.1}, then the conclusion also holds for $N = n$.

\begin{theorem} 
\label{higher.dim.inequality.2}
Let us fix an integer $n$ with $3 \leq n \leq 7$. Let $M$ be a compact, connected, orientable manifold of dimension $n$ with non-empty boundary $\partial M$. Suppose that $\theta_0$ is a smooth map from $\partial M$ to $S^1$, and $(\theta_1,\hdots,\theta_{n-2})$ is a smooth map from $M$ to $T^{n-2}$. We assume that the map $(\theta_0,\theta_1|_{\partial M},\hdots,\theta_{n-2}|_{\partial M}): \partial M \to S^1 \times T^{n-2}$ has non-zero degree. We denote by $\Theta_0$ the pull-back of the volume form on $S^1$ under the map $\theta_0: \partial M \to S^1$. Note that $\Theta_0$ is a closed one-form on $\partial M$. Let $g$ be a Riemannian metric on $M$ with scalar curvature at least $-n(n-1)$, and let $H_{\partial M}$ denote the mean curvature of the boundary $\partial M$. Then 
\[2 \, \sigma^n \, \inf_{\partial M} (H_{\partial M} - (n-1)) < \Big ( \frac{4\pi}{n} \Big )^n,\] 
where $\sigma$ denotes the length of the shortest closed curve $\alpha$ in $\partial M$ satisfying $\int_\alpha \Theta_0 \neq 0$.
\end{theorem}

Theorem \ref{higher.dim.inequality.1} can be compared to the following theorem (see also \cite{Ambrozio}).

\begin{theorem} 
\label{higher.dim.inequality.3}
Let us fix an integer $n$ with $3 \leq n \leq 7$. Let $M$ be a compact, connected, orientable manifold of dimension $n$ with non-empty boundary $\partial M$. Suppose that $\theta_0$ is a smooth map from $\partial M$ to $S^1$, and $(\theta_1,\hdots,\theta_{n-2})$ is a smooth map from $M$ to $T^{n-2}$. We assume that the map $(\theta_0,\theta_1|_{\partial M},\hdots,\theta_{n-2}|_{\partial M}): \partial M \to S^1 \times T^{n-2}$ has non-zero degree. We denote by $\Theta_0$ the pull-back of the volume form on $S^1$ under the map $\theta_0: \partial M \to S^1$. Note that $\Theta_0$ is a closed one-form on $\partial M$. Let $g$ be a Riemannian metric on $M$. We denote by $R_M$ the scalar curvature of $M$, by $H_{\partial M}$ the mean curvature of the boundary $\partial M$, and by $\eta$ the outward-pointing unit normal vector field to $\partial M$. Let $\varphi$ be a smooth function on $M$ such that 
\[-2 \, \Delta_M \varphi - |\nabla^M \varphi|^2 + R_M \geq 0\] 
at each point in $M$. Then 
\[\sigma \, \inf_{\partial M} (\langle \nabla^M \varphi,\eta \rangle + H_{\partial M}) \leq 2\pi,\] 
where $\sigma$ denotes the length of the shortest closed curve $\alpha$ in $\partial M$ satisfying $\int_\alpha \Theta_0 \neq 0$.
\end{theorem}

The case of equality in Theorem \ref{higher.dim.inequality.3} is studied in \cite{Wang2}.

As a consequence of Theorem \ref{higher.dim.inequality.2}, we obtain a new positive energy theorem, as envisioned by Horowitz and Myers \cite{Horowitz-Myers}. 

\begin{theorem} 
\label{positive.energy.theorem}
Let us fix an integer $n$ with $3 \leq n \leq 7$. Let $\Theta_0$ denote the pull-back of the volume form on $S^1$ under the canonical projection from $S^1 \times T^{n-2}$ to $S^1$. Note that $\Theta_0$ is a closed one-form on $S^1 \times T^{n-2}$. Let us fix a flat metric $\gamma$ on $S^1 \times T^{n-2}$. Let $Q$ be a symmetric $(0,2)$-tensor on $S^1 \times T^{n-2}$. Given a positive real number $r_0$, we define a hyperbolic metric $\bar{g}$ on $[r_0,\infty) \times S^1 \times T^{n-2}$ by $\bar{g} = r^{-2} \, dr \otimes dr + r^2 \, \gamma$. Let $(M,g)$ be a complete, connected, orientable Riemannian manifold of dimension $n$ with the following properties: 
\begin{itemize}
\item There exists a bounded open domain $E \subset M$ with smooth boundary such that the complement $M \setminus E$ is diffeomorphic to $[r_0,\infty) \times S^1 \times T^{n-2}$.
\item The canonical projection from $[r_0,\infty) \times S^1 \times T^{n-2}$ to $T^{n-2}$ extends to a smooth map from $M$ to $T^{n-2}$.
\item On the complement $M \setminus E$, the metric satisfies 
\[|g-\bar{g} - r^{2-n} \, Q|_{\bar{g}} \leq o(r^{-n})\] 
and 
\[|\bar{D}(g-\bar{g} - r^{2-n} \, Q)|_{\bar{g}} \leq o(r^{-n})\] 
as $r \to \infty$. 
\item The metric $g$ has scalar curvature at least $-n(n-1)$ at each point in $M$. 
\end{itemize}
Then 
\[\int_{S^1 \times T^{n-2}} \Big ( n \, \text{\rm tr}_\gamma(Q) + \Big ( \frac{4\pi}{n\sigma} \Big )^n \Big ) \, d\text{\rm vol}_\gamma \geq 0,\] 
where $\sigma$ denotes the length of the shortest closed curve $\alpha$ in $(S^1 \times T^{n-2},\gamma)$ satisfying $\int_\alpha \Theta_0 \neq 0$.
\end{theorem}

We refer to \cite{Barzegar-Chrusciel-Horzinger-Maliborski-Nguyen}, \cite{Chrusciel-Delay1}, \cite{Chrusciel-Delay-Wutte}, \cite{Chrusciel-Galloway-Nguyen-Paetz}, \cite{Constable-Myers}, \cite{Galloway-Woolgar}, \cite{Horowitz-Myers}, \cite{Liang-Zhang}, \cite{Woolgar} for previous work on the Horowitz-Myers conjecture, and to \cite{Andersson-Cai-Galloway}, \cite{Andersson-Dahl}, \cite{Chrusciel-Delay2}, \cite{Chrusciel-Galloway}, \cite{Lee-Neves}, \cite{Min-Oo}, \cite{Wang1} for other scalar curvature rigidity results for asymptotically hyperbolic manifolds. In a remarkable work, Barzegar, Chru\'sciel, H\"orzinger, Maliborski, and Nguyen \cite{Barzegar-Chrusciel-Horzinger-Maliborski-Nguyen} verified the Horowitz-Myers conjecture in the special case when $(M,g)$ is a multiply warped product over a halfline.

\begin{remark} 
The proof of Theorem \ref{positive.energy.theorem} does not use the assumption that $\gamma$ is flat. If $\gamma$ is an arbitrary metric on $S^1 \times T^{n-2}$, then the metric $\bar{g} = r^{-2} \, dr \otimes dr + r^2 \, \gamma$ is asymptotically hyperbolic, and the scalar curvature of $\bar{g}$ is related to the scalar curvature of $\gamma$ by the formula $R_{\bar{g}} = -n(n-1) + r^{-2} \, R_\gamma$. Hence, if $R_{\bar{g}} \geq -n(n-1) - o(r^{-2})$, then $R_\gamma \geq 0$, and this implies that $\gamma$ is flat by the solution of Geroch's conjecture (see \cite{Gromov-Lawson}, \cite{Schoen-Yau1}). 
\end{remark}

\begin{remark}
The inequality in Theorem \ref{positive.energy.theorem} is sharp for the Horowitz-Myers metrics. To see this, let us fix an integer $n \geq 3$ and a flat metric $g_{T^{n-2}}$ on $T^{n-2}$. Let $g_{S^1}$ denote the standard metric on $S^1$, so that $(S^1,g_{S^1})$ has length $2\pi$. We define a metric $g$ on $(2^{-\frac{2}{n}},\infty) \times S^1 \times T^{n-2}$ by 
\begin{align*} 
g &= r^{-2} \, dr \otimes dr + \frac{4}{n^2} \, r^2 \, (1 + \frac{1}{4} \, r^{-n})^{\frac{4}{n}-2} \, (1 - \frac{1}{4} \, r^{-n})^2 \, g_{S^1} \\ 
&+ r^2 \, (1 + \frac{1}{4} \, r^{-n})^{\frac{4}{n}} \, g_{T^{n-2}}, 
\end{align*} 
where $r \in (2^{-\frac{2}{n}},\infty)$. The metric $g$ extends to a smooth metric on $\mathbb{R}^2 \times T^{n-2}$. We define a metric $\gamma$ on $S^1 \times T^{n-2}$ by 
\[\gamma = \frac{4}{n^2} \, g_{S^1} + g_{T^{n-2}}.\] 
Moreover, we define a symmetric $(0,2)$-tensor $Q$ on $S^1 \times T^{n-2}$ by 
\[Q = -\frac{4(n-1)}{n^3} \, g_{S^1} + \frac{1}{n} \, g_{T^{n-2}}.\] 
As above, we consider the hyperbolic metric $\bar{g} = r^{-2} \, dr \otimes dr + r^2 \, \gamma$. From the definition of $g$, it is easy to see that 
\[|g-\bar{g}-r^{2-n} \, Q|_{\bar{g}} \leq o(r^{-n})\] 
as $r \to \infty$. Moreover, $\text{\rm tr}_\gamma(Q) = -\frac{1}{n}$ at each point on $S^1 \times T^{n-2}$. On the other hand, since $(S^1,g_{S^1})$ has length $2\pi$, it follows that $\sigma = \frac{4\pi}{n}$. This gives 
\[n \, \text{\rm tr}_\gamma(Q) + \Big ( \frac{4\pi}{n\sigma} \Big )^n = 0.\] 
Thus, the metric $g$ achieves equality in Theorem \ref{positive.energy.theorem}.
\end{remark}

\begin{remark} 
Using the substitution $\Upsilon = r \, (1 + \frac{1}{4} \, r^{-n})^{\frac{2}{n}}$, the Horowitz-Myers metric can be written in the form 
\[g = \Upsilon^{-2} \, (1-\Upsilon^{-n})^{-1} \, d\Upsilon \otimes d\Upsilon + \frac{4}{n^2} \, \Upsilon^2 \, (1-\Upsilon^{-n}) \, g_{S^1} + \Upsilon^2 \, g_{T^{n-2}},\] 
with $\Upsilon$ taking values in the interval $[1,\infty)$. A straightforward calculation shows that $g$ is a static metric with scalar curvature $-n(n-1)$ (see \cite{Woolgar}). 
\end{remark}

\section{Proof of Theorem \ref{2D.inequality} -- The case when $\partial \Sigma$ is connected}

In this section, we prove Theorem \ref{2D.inequality} in the special case when $\partial \Sigma$ is connected. Let $N > 2$ be a real number, let $\Sigma$ be a two-dimensional surface with non-empty boundary $\partial \Sigma$, let $g$ be a Riemannian metric on $\Sigma$, and let $\psi$ be a smooth function on $\Sigma$. We assume that $\Sigma$, $g$, and $\psi$ satisfy the assumptions of Theorem \ref{2D.inequality} and that the boundary $\partial \Sigma$ is connected. We define a function $w: \Sigma \to [0,\infty)$ by 
\[w(x) = d(x,\partial \Sigma)\] 
for each point $x \in \Sigma$. Let 
\[l = \sup_{y \in \Sigma} w(y).\] 
As in Section 4.2 of \cite{Shiohama-Shioya-Tanaka}, we define a function $\rho: \partial \Sigma \to (0,\infty)$ by 
\[\rho(x) = \sup \big \{ s>0: \text{\rm $\exp_x(-s \, \eta(x))$ is defined and $w(\exp_x(-s \, \eta(x))) = s$} \big \}\] 
for each point $x \in \partial \Sigma$. Clearly, $\rho(x) \leq l$ for all $x \in \partial \Sigma$. Let 
\[\text{\rm seg} = \{(x,s) \in \partial \Sigma \times [0,\infty): s \in [0,\rho(x)]\}\] 
and 
\[\text{\rm seg}^0 = \{(x,s) \in \partial \Sigma \times [0,\infty): s \in [0,\rho(x))\}.\] 
It is well known that $\text{\rm seg}$ is a compact subset of $\partial \Sigma \times [0,\infty)$ and $\text{\rm seg}^0$ is a relatively open subset of $\partial \Sigma \times [0,\infty)$ (compare \cite{Petersen}, Proposition 5.7.10). Thus, the function $\rho: \partial \Sigma \to (0,\infty)$ is continuous (see also \cite{Shiohama-Shioya-Tanaka}, Proposition 4.2.1). 

We define a map $\Phi: \text{\rm seg} \to \Sigma$ by 
\[\Phi(x,s) = \exp_x(-s \, \eta(x))\] 
for all pairs $(x,s) \in \text{\rm seg}$. Note that $\Phi$ is surjective. Moreover, if $(x,s) \in \text{\rm seg}^0$ and $(\tilde{x},\tilde{s}) \in \text{\rm seg}$ satisfy $\Phi(x,s) = \Phi(\tilde{x},\tilde{s})$, then $(x,s) = (\tilde{x},\tilde{s})$ (compare \cite{Petersen}, Proposition 5.7.7). In particular, the restriction $\Phi|_{\text{\rm seg}^0}$ is injective. Finally, with a suitable choice of orientation, we have $\det (D\Phi)_{(x,s)} > 0$ for each point $(x,s) \in \text{\rm seg}^0$.

For each $s \in (0,l)$, we denote by $A(s)$ the area of the tubular neighborhood $\{y \in \Sigma: w(y) \leq s\}$. We may write 
\[A(s) = \int_{\partial \Sigma} \bigg ( \int_0^{\min \{\rho(x),s\}} \det (D\Phi)_{(x,t)} \, dt \bigg ) \, d\mathcal{H}^1(x)\] 
for each $s \in (0,l)$. Moreover, we define 
\[L(s) = \int_{\partial \Sigma} 1_{\{\rho(x) \geq s\}} \det (D\Phi)_{(x,s)} \, d\mathcal{H}^1(x)\] 
for each $s \in (0,l)$. Note that the function $s \mapsto L(s)$ is not necessarily continuous.

\begin{lemma} 
\label{A.almost.monotone.decreasing}
We can find a large constant $C_1$ such that the function $s \mapsto A(s) - C_1 \, s$ is monotone decreasing for $s \in (0,l)$. Moreover, $|\Sigma| - A(s) \leq C_1 \, (l-s)$ for each $s \in (0,l)$. 
\end{lemma} 

\textbf{Proof.} 
Let us fix a large constant $C_2$ such that $\det (D\Phi)_{(x,s)} \leq C_2$ for each point $x \in \partial \Sigma$ and each $s \in [0,\rho(x)]$. This implies 
\begin{align*}
A(s_1) - A(s_0) 
&= \int_{\partial \Sigma} \bigg ( \int_{\min \{\rho(x),s_0\}}^{\min \{\rho(x),s_1\}} \det (D\Phi)_{(x,t)} \, dt \bigg ) \, d\mathcal{H}^1(x) \\ 
&\leq C_2 \, |\partial \Sigma| \, (s_1-s_0) 
\end{align*} 
for $0 < s_0 < s_1 < l$. Moreover, using the identity 
\[|\Sigma| = \int_{\partial \Sigma} \bigg ( \int_0^{\rho(x)} \det (D\Phi)_{(x,t)} \, dt \bigg ) \, d\mathcal{H}^1(x)\] 
and the fact that $\rho(x) \leq l$ for all $x \in \partial \Sigma$, we obtain 
\begin{align*}
|\Sigma| - A(s) 
&= \int_{\partial \Sigma} \bigg ( \int_{\min \{\rho(x),s\}}^{\rho(x)} \det (D\Phi)_{(x,t)} \, dt \bigg ) \, d\mathcal{H}^1(x) \\ 
&\leq C_2 \, |\partial \Sigma| \, (l-s) 
\end{align*} 
for each $s \in (0,l)$. This completes the proof of Lemma \ref{A.almost.monotone.decreasing}. \\

\begin{lemma} 
\label{L.almost.monotone.decreasing}
We can find a large constant $C_3$ such that the function $s \mapsto L(s) - C_3 \, s$ is monotone decreasing for $s \in (0,l)$.
\end{lemma} 

\textbf{Proof.} 
Standard results in comparison geometry imply that we can find a large constant $C_4$ such that 
\[\frac{\partial}{\partial s} \log \det (D\Phi)_{(x,s)} \leq C_4\] 
for each point $x \in \partial \Sigma$ and each $s \in [0,\rho(x))$. For each point $x \in \partial \Sigma$, the function $s \mapsto e^{-C_4 s} \, \det (D\Phi)_{(x,s)}$ is monotone decreasing for $s \in [0,\rho(x)]$. Consequently, the function $s \mapsto e^{-C_4 s} \, L(s)$ is monotone decreasing for $s \in (0,l)$. From this, the assertion follows easily. This completes the proof of Lemma \ref{L.almost.monotone.decreasing}. \\

We define a function $F: (0,\infty) \to (0,1)$ by 
\[F(s) = \tanh \Big ( \frac{Ns}{2} \Big )\] 
for each $s \in (0,\infty)$. Moreover, we define a function $G: (0,\infty) \to (1,\infty)$ by 
\[G(s) = \Big [ \cosh \Big ( \frac{Ns}{2} \Big ) \Big ]^{\frac{2(N-1)}{N}}\] 
for each $s \in (0,\infty)$. The derivatives of $F$ and $G$ are given by 
\begin{equation} 
\label{derivative.F}
F'(s) = \frac{N}{2} \, (1-F(s)^2)
\end{equation} 
and 
\begin{equation} 
\label{derivative.G}
G'(s) = (N-1) \, G(s) \, F(s)
\end{equation} 
for each $s \in (0,\infty)$. 

For each $s \in (0,l)$, we define 
\[\Omega_s = \{y \in \Sigma: w(y) > s\}.\] 
For each $s \in (0,l)$, $\Omega_s$ is a non-empty open subset of $\Sigma$ with area $|\Sigma| - A(s)$. Finally, we define 
\[I(s) = 2\pi \chi(\Sigma) - (N-1) \, F(l-s) \, L(s) + \int_{\Omega_s} (\Delta \psi - K)\] 
and 
\[J(s) = G(l-s) \, I(s)\] 
for each $s \in (0,l)$. 

\begin{lemma}
\label{I.almost.monotone.increasing}
We can find a large constant $C_5$ such that $I(s) \leq C_5$ for all $s \in (0,l)$. Moreover, we can find a large constant $C_6$ such that the function $s \mapsto I(s) + C_6 \, s$ is monotone increasing for $s \in (0,l)$.
\end{lemma} 

\textbf{Proof.} 
We may bound 
\[I(s) \leq 2\pi \, \chi(\Sigma) + \int_\Sigma |\Delta \psi - K|\] 
for all $s \in (0,l)$. This proves the first statement.

It follows from Lemma \ref{A.almost.monotone.decreasing} and Lemma \ref{L.almost.monotone.decreasing} that the functions $s \mapsto A(s) - C_1 \, s$ and $s \mapsto L(s) - C_3 \, s$ are monotone decreasing for $s \in (0,l)$. This implies 
\[L(s_1)-L(s_0) \leq C_3 \, (s_1-s_0)\] 
and 
\begin{align*} 
\int_{\Omega_{s_0} \setminus \Omega_{s_1}} (\Delta \psi - K) 
&\leq \sup_\Sigma |\Delta \psi - K| \, (A(s_1)-A(s_0)) \\ 
&\leq C_1 \, \sup_\Sigma |\Delta \psi - K| \, (s_1-s_0) 
\end{align*}
for $0 < s_0 < s_1 < l$. Since the function $F$ is monotone increasing, we obtain 
\[0 \leq F(l-s_1) \leq F(l-s_0) \leq 1\] 
for $0 < s_0 < s_1 < l$. Putting these facts together, we conclude that 
\begin{align*} 
I(s_1) - I(s_0) 
&= -(N-1) \, F(l-s_0) \, (L(s_1) - L(s_0)) \\ 
&+ (N-1) \, (F(l-s_0) - F(l-s_1)) \, L(s_1) \\ 
&- \int_{\Omega_{s_0} \setminus \Omega_{s_1}} (\Delta \psi - K) \\ 
&\geq -(N-1) \, C_3 \, (s_1-s_0) \\ 
&- C_1 \, \sup_\Sigma |\Delta \psi - K| \, (s_1-s_0) 
\end{align*}
for $0 < s_0 < s_1 < l$. This proves the second statement. This completes the proof of Lemma \ref{I.almost.monotone.increasing}. \\

\begin{lemma}
\label{J.almost.monotone.increasing}
We can find a large constant $C_7$ such that $J(s) \leq C_7$ for all $s \in (0,l)$. Moreover, we can find a large constant $C_8$ such that the function $s \mapsto J(s) + C_8 \, s$ is monotone increasing for $s \in (0,l)$.
\end{lemma}

\textbf{Proof.} 
It follows from Lemma \ref{I.almost.monotone.increasing} that $I(s) \leq C_5$ for all $s \in (0,l)$. This implies 
\[J(s) = G(l-s) \, I(s) \leq C_5 \, G(l-s) \leq C_5 \, G(l)\] 
for all $s \in (0,l)$. This proves the first statement.

By Lemma \ref{I.almost.monotone.increasing}, the function $s \mapsto I(s) + C_6 \, s$ is monotone increasing for $s \in (0,l)$. This gives 
\[I(s_1) - I(s_0) \geq -C_6 \, (s_1-s_0)\] 
for $0 < s_0 < s_1 < l$. Using (\ref{derivative.G}), we obtain $0 \leq G' \leq (N-1) \, G(l)$ on the interval $(0,l)$. Consequently, 
\[0 \leq G(l-s_0) - G(l-s_1) \leq (N-1) \, G(l) \, (s_1-s_0)\] 
for $0 < s_0 < s_1 < l$. Putting these facts together, we conclude that 
\begin{align*} 
J(s_1)-J(s_0) 
&= G(l-s_0) \, (I(s_1)-I(s_0)) - (G(l-s_0)-G(l-s_1)) \, I(s_1) \\ 
&\geq -C_6 \, G(l-s_0) \, (s_1-s_0) - C_5 \, (G(l-s_0)-G(l-s_1)) \\ 
&\geq -(C_6 + (N-1) \, C_5) \, G(l) \, (s_1-s_0) 
\end{align*} 
for $0 < s_0 < s_1 < l$. This proves the second statement. This completes the proof of Lemma \ref{J.almost.monotone.increasing}. \\

The following result was proved by Shiohama, Shioya, and Tanaka \cite{Shiohama-Shioya-Tanaka}; see \cite{Fiala} and \cite{Hartman} for related work.

\begin{theorem}[cf. K.~Shiohama, T.~Shioya, M.~Tanaka \cite{Shiohama-Shioya-Tanaka}]
\label{generic.level.sets}
We can find a relatively closed set $\mathcal{E} \subset (0,l)$ of measure zero with the following properties: 
\begin{itemize}
\item[(i)] Suppose that $s \in (0,l) \setminus \mathcal{E}$. Moreover, suppose that $x$ is a point in $\partial \Sigma$ with $\rho(x) \geq s$. Then $\det (D\Phi)_{(x,s)} \neq 0$. 
\item[(ii)] Suppose that $s \in (0,l) \setminus \mathcal{E}$. Moreover, suppose that $y$ is a point in $\Sigma$ with $w(y) = s$. Then the set $\{x \in \partial \Sigma: \text{\rm $\rho(x) \geq s$ and $\Phi(x,s) = y$}\}$ consists of at most two elements. 
\item[(iii)] Suppose that $s \in (0,l) \setminus \mathcal{E}$. Then the set $\{x \in \partial \Sigma: \rho(x) = s\}$ is finite.
\item[(iv)] Suppose that $s \in (0,l) \setminus \mathcal{E}$. Then $\Omega_s$ is a domain with piecewise smooth boundary. The length of the boundary $\partial \Omega_s$ is equal to $L(s)$. Moreover, 
\begin{equation} 
\label{derivative.of.L}
\limsup_{\delta \searrow 0} \frac{L(s+\delta)-L(s)}{\delta} \leq -\Lambda(s), 
\end{equation}
where $\Lambda(s)$ denotes the total geodesic curvature of the boundary $\partial \Omega_s$ (including angle contributions).
\end{itemize}
\end{theorem} 

\textbf{Proof.} 
These statements follow from results in Chapter 4 of \cite{Shiohama-Shioya-Tanaka}. Let $\mathcal{E}$ denote the set of all real numbers $s \in (0,l)$ with the property that $s$ is an exceptional value in the sense of Definition 4.3.1 in \cite{Shiohama-Shioya-Tanaka}. By Lemma 4.3.6 in \cite{Shiohama-Shioya-Tanaka}, the set $\mathcal{E}$ is relatively closed and has measure zero. 

We claim that the set $\mathcal{E}$ satisfies the conditions (i)--(iv) above. To prove this, we fix a real number $s \in (0,l) \setminus \mathcal{E}$. In other words, $s$ is non-exceptional in the sense of Definition 4.3.1 in \cite{Shiohama-Shioya-Tanaka}. In particular, $s$ is normal in the sense of Definition 4.3.1 in \cite{Shiohama-Shioya-Tanaka} (see also p.~142 in \cite{Shiohama-Shioya-Tanaka}). From this, properties (i) and (ii) follow immediately. Property (iii) follows from Lemma 4.4.1 in \cite{Shiohama-Shioya-Tanaka}. It remains to prove property (iv). By Theorem 4.4.1 in \cite{Shiohama-Shioya-Tanaka}, $\Omega_s$ is a domain with piecewise smooth boundary. It follows from equation (4.4.1) in \cite{Shiohama-Shioya-Tanaka} that the length of the boundary $\partial \Omega_s$ is equal to $L(s)$. Moreover, the interior angles of the domain $\Omega_s$ are at most $\pi$, and the angles $\theta_k$ defined in Theorem 4.4.1 in \cite{Shiohama-Shioya-Tanaka} can be interpreted as the exterior angles of $\Omega_s$. Note that $\theta_k \in [0,\pi)$ for each $k$. This implies $\tan \frac{\theta_k}{2} \geq \frac{\theta_k}{2}$ for each $k$. Hence, the inequality (\ref{derivative.of.L}) follows from equation (4.4.2) in \cite{Shiohama-Shioya-Tanaka}. This completes the proof of Theorem \ref{generic.level.sets}. \\

\begin{lemma}
\label{Euler.characteristic}
For each $s \in (0,l) \setminus \mathcal{E}$, the Euler characteristic of the domain $\Omega_s$ satisfies $\chi(\Omega_s) \geq \chi(\Sigma)$.
\end{lemma} 

\textbf{Proof.} 
Let us fix a real number $s \in (0,l) \setminus \mathcal{E}$. Since the boundary $\partial \Sigma$ is connected, it follows that the tubular neighborhood 
\[\Sigma \setminus \Omega_s = \{y \in \Sigma: w(y) \leq s\}\] 
is connected as well. Thus, $\Sigma \setminus \Omega_s$ is homeomorphic to a compact, connected, orientable manifold of dimension $2$ with disconnected boundary. Consequently, the Euler characteristic of $\Sigma \setminus \Omega_s$ is non-positive. Since $\chi(\Sigma \setminus \Omega_s) + \chi(\Omega_s) = \chi(\Sigma)$, we conclude that $\chi(\Omega_s) \geq \chi(\Sigma)$. This completes the proof of Lemma \ref{Euler.characteristic}. \\

\begin{lemma} 
\label{auxiliary.result}
For each $s \in (0,l) \setminus \mathcal{E}$, we have 
\[2\pi \chi(\Sigma) - \int_{\Omega_s} K \leq \Lambda(s)\] 
and 
\[\int_{\Omega_s} \Delta \psi \leq \int_{\partial \Omega_s} |\nabla \psi|.\] 
As above, $\Lambda(s)$ denotes the total geodesic curvature of the boundary $\partial \Omega_s$ (including angle contributions).
\end{lemma} 

\textbf{Proof.} 
Let us fix a real number $s \in (0,l) \setminus \mathcal{E}$. It follows from Lemma \ref{Euler.characteristic} that $\chi(\Omega_s) \geq \chi(\Sigma)$. Hence, the first statement follows from the Gauss-Bonnet theorem. The second statement follows from the divergence theorem. This completes the proof of Lemma \ref{auxiliary.result}. \\

\begin{proposition}
\label{derivative.of.I}
For each $s \in (0,l) \setminus \mathcal{E}$, we have 
\[\liminf_{\delta \searrow 0} \frac{I(s+\delta)-I(s)}{\delta} - (N-1) \, F(l-s) \, I(s) \geq 0.\] 
\end{proposition} 

\textbf{Proof.} 
Let us fix a real number $s \in (0,l) \setminus \mathcal{E}$. By Theorem \ref{generic.level.sets} (iv), we obtain 
\[\limsup_{\delta \searrow 0} \frac{L(s+\delta)-L(s)}{\delta} \leq -\Lambda(s).\]
Using (\ref{derivative.F}), it follows that 
\begin{align*} 
\liminf_{\delta \searrow 0} \frac{I(s+\delta)-I(s)}{\delta} 
&\geq \frac{N(N-1)}{2} \, (1-F(l-s)^2) \, L(s) \\ 
&+ (N-1) \, F(l-s) \, \Lambda(s) - \int_{\partial \Omega_s} (\Delta \psi - K). 
\end{align*} 
On the other hand, Lemma \ref{auxiliary.result} gives 
\[I(s) \leq -(N-1) \, F(l-s) \, L(s) + \Lambda(s) + \int_{\partial \Omega_s} |\nabla \psi|.\] 
Putting these facts together, we conclude that 
\begin{align*} 
&\liminf_{\delta \searrow 0} \frac{I(s+\delta)-I(s)}{\delta} - (N-1) \, F(l-s) \, I(s) \\ 
&\geq \frac{N(N-1)}{2} \, L(s) + \frac{(N-2)(N-1)}{2} \, F(l-s)^2 \, L(s) \\ 
&- (N-1) \, F(l-s) \int_{\partial \Omega_s} |\nabla \psi| - \int_{\partial \Omega_s} (\Delta \psi - K). 
\end{align*} 
Rearranging terms gives 
\begin{align*} 
&\liminf_{\delta \searrow 0} \frac{I(s+\delta)-I(s)}{\delta} - (N-1) \, F(l-s) \, I(s) \\ 
&\geq \frac{N-1}{2(N-2)} \int_{\partial \Omega_s} ((N-2) \, F(l-s) - |\nabla \psi|)^2 \\ 
&+ \int_{\partial \Omega_s} \Big ( -\Delta \psi - \frac{N-1}{2(N-2)} \, |\nabla \psi|^2 + K + \frac{N(N-1)}{2} \Big ). 
\end{align*} 
The assertion therefore follows from the assumption that 
\[-\Delta \psi - \frac{N-1}{2(N-2)} \, |\nabla \psi|^2 + K + \frac{N(N-1)}{2} \geq 0\] 
at each point in $\Sigma$. This completes the proof of Proposition \ref{derivative.of.I}. \\

\begin{proposition}
\label{derivative.of.J}
For each $s \in (0,l) \setminus \mathcal{E}$, we have 
\[\liminf_{\delta \searrow 0} \frac{J(s+\delta)-J(s)}{\delta} \geq 0.\] 
\end{proposition} 

\textbf{Proof.} 
Let us fix a real number $s \in (0,l) \setminus \mathcal{E}$. Using (\ref{derivative.G}), we obtain 
\begin{align*} 
\liminf_{\delta \searrow 0} \frac{J(s+\delta)-J(s)}{\delta} 
&= G(l-s) \liminf_{\delta \searrow 0} \frac{I(s+\delta)-I(s)}{\delta} \\ 
&- (N-1) \, G(l-s) \, F(l-s) \, I(s), 
\end{align*}
and the expression on the right hand side is nonnegative by Proposition \ref{derivative.of.I}. This completes the proof of Proposition \ref{derivative.of.J}. \\

\begin{corollary} 
\label{monotonicity} 
The function $s \mapsto J(s)$ is monotone increasing for $s \in (0,l)$.
\end{corollary} 

\textbf{Proof.} 
It follows from Lemma \ref{J.almost.monotone.increasing} and the monotone differentiation theorem (see e.g. \cite{Royden}, p.~100) that the function $s \mapsto J(s)$ is differentiable almost everywhere and 
\[J(s_1) - J(s_0) \geq \int_{s_0}^{s_1} J'(s) \, ds\] 
for all $0 < s_0 < s_1 < l$. Since $\mathcal{E}$ is a set of measure zero, Proposition \ref{derivative.of.J} implies that $J'(s) \geq 0$ almost everywhere. Putting these facts together, we conclude that $J(s_1) - J(s_0) \geq 0$ for all $0 < s_0 < s_1 < l$. This completes the proof of Corollary \ref{monotonicity}. \\

\begin{proposition} 
\label{consequence.of.monotonicity.1}
We have 
\[\int_{\partial \Sigma} (\langle \nabla \psi,\eta \rangle + \kappa) \leq 2\pi \chi(\Sigma) \, G(l)^{-1} + (N-1) \, F(l) \, |\partial \Sigma|.\] 
\end{proposition}

\textbf{Proof.}  
Using the Gauss-Bonnet theorem and the divergence theorem, we obtain 
\begin{align*} 
\liminf_{s \searrow 0} I(s) 
&= 2\pi \chi(\Sigma) - (N-1) \, F(l) \, |\partial \Sigma| + \int_\Sigma (\Delta \psi - K) \\ 
&= -(N-1) \, F(l) \, |\partial \Sigma| + \int_{\partial \Sigma} (\langle \nabla \psi,\eta \rangle + \kappa). 
\end{align*} 
This gives 
\[\liminf_{s \searrow 0} J(s) = G(l) \, \bigg ( -(N-1) \, F(l) \, |\partial \Sigma| + \int_{\partial \Sigma} (\langle \nabla \psi,\eta \rangle + \kappa) \bigg ).\] 
On the other hand, 
\[I(s) \leq 2\pi \chi(\Sigma) + \int_{\Omega_s} |\Delta \psi - K| \leq 2\pi \chi(\Sigma) + \sup_\Sigma |\Delta \psi - K| \, (|\Sigma| - A(s))\] 
for each $s \in (0,l)$. Moreover, Lemma \ref{A.almost.monotone.decreasing} implies that $|\Sigma| - A(s) \leq C_1 \, (l-s)$ for each $s \in (0,l)$. Thus, 
\[\limsup_{s \nearrow l} I(s) \leq 2\pi \chi(\Sigma).\] 
Since $\lim_{s \nearrow l} G(l-s) = 1$, it follows that 
\[\limsup_{s \nearrow l} J(s) \leq 2\pi \chi(\Sigma).\] 
Using Corollary \ref{monotonicity}, we conclude that 
\[G(l) \, \bigg ( -(N-1) \, F(l) \, |\partial \Sigma| + \int_{\partial \Sigma} (\langle \nabla \psi,\eta \rangle + \kappa) \bigg ) \leq 2\pi \chi(\Sigma).\] 
This completes the proof of Proposition \ref{consequence.of.monotonicity.1}. \\

\begin{proposition}
\label{consequence.of.monotonicity.2}
We have 
\[2 \, |\partial \Sigma|^{N-1} \, \int_{\partial \Sigma} (\langle \nabla \psi,\eta \rangle + \kappa - (N-1)) \leq \Big ( \frac{4\pi}{N} \Big )^N - (N-1) \, (1-F(l))^2 \, |\partial \Sigma|^N.\] 
\end{proposition} 

\textbf{Proof.} 
By assumption, $\Sigma$ is a compact, connected, orientable manifold of dimension $2$ with non-empty boundary. Consequently, $\chi(\Sigma) \leq 1$. Using Proposition \ref{consequence.of.monotonicity.1} and the identity $1-F(l)^2 = G(l)^{-\frac{N}{N-1}}$, we obtain 
\begin{align*} 
&2 \, |\partial \Sigma|^{N-1} \, \int_{\partial \Sigma} (\langle \nabla \psi,\eta \rangle + \kappa - (N-1)) \\ 
&\leq 4\pi \, G(l)^{-1} \, |\partial \Sigma|^{N-1} - 2(N-1) \, (1-F(l)) \, |\partial \Sigma|^N \\ 
&= 4\pi \, G(l)^{-1} \, |\partial \Sigma|^{N-1} - (N-1) \, G(l)^{-\frac{N}{N-1}} \, |\partial \Sigma|^N \\ 
&- (N-1) \, (1-F(l))^2 \, |\partial \Sigma|^N. 
\end{align*} 
It follows from Young's inequality that $Nab \leq a^N + (N-1) \, b^{\frac{N}{N-1}}$ for all real numbers $a,b \geq 0$. Putting $a = \frac{4\pi}{N}$ and $b = G(l)^{-1} \, |\partial \Sigma|^{N-1}$ gives 
\[4\pi \, G(l)^{-1} \, |\partial \Sigma|^{N-1} \leq \Big ( \frac{4\pi}{N} \Big )^N + (N-1) \, G(l)^{-\frac{N}{N-1}} \, |\partial \Sigma|^N.\] 
Putting these facts together, the assertion follows. This completes the proof of Proposition \ref{consequence.of.monotonicity.2}. \\

Proposition \ref{consequence.of.monotonicity.2} implies 
\[2 \, |\partial \Sigma|^N \, \inf_{\partial \Sigma} (\langle \nabla \psi,\eta \rangle + \kappa - (N-1)) \leq \Big ( \frac{4\pi}{N} \Big )^N - (N-1) \, (1-F(l))^2 \, |\partial \Sigma|^N.\] 
This completes the proof of Theorem \ref{2D.inequality} in the special case when $\partial \Sigma$ is connected.

\section{Proof of Theorem \ref{2D.inequality} -- The case when $\partial \Sigma$ is disconnected}

In this section, we prove Theorem \ref{2D.inequality} in the special case when $\partial \Sigma$ is disconnected. Let $N > 2$ be a real number, let $\Sigma$ be a two-dimensional surface with non-empty boundary $\partial \Sigma$, let $g$ be a Riemannian metric on $\Sigma$, and let $\psi$ be a smooth function on $\Sigma$. We assume that $\Sigma$, $g$, and $\psi$ satisfy the assumptions of Theorem \ref{2D.inequality} and that the boundary $\partial \Sigma$ is disconnected. We assume that the conclusion of Theorem \ref{2D.inequality} is false, so that 
\begin{equation} 
\label{lower.bound.for.weighted.geodesic.curvature}
\inf_{\partial \Sigma} (\langle \nabla \psi,\eta \rangle + \kappa) > N-1. 
\end{equation}
In particular, the surface $(\Sigma,e^{2\psi} \, g)$ has strictly convex boundary. 

We minimize the weighted energy functional $\int_0^1 e^{2\psi(\beta(t))} \, |\beta'(t)|^2 \, dt$ among all smooth maps $\beta: [0,1] \to \Sigma$ with the property that $\beta(0),\beta(1) \in \partial \Sigma$ and $\beta(0)$ and $\beta(1)$ lie in different connected components of $\partial \Sigma$. It is well known that the minimum is attained by a smooth map $\beta: [0,1] \to \Sigma$. Moreover, $\beta$ is a geodesic in $(\Sigma,e^{2\psi} \, g)$ which meets the boundary $\partial \Sigma$ orthogonally. 

Let $\gamma: [0,l] \to \Sigma$ be a reparametrization of $\beta$ with the property that $\gamma$ has unit speed with respect to the metric $g$. Note that $\gamma$ minimizes the weighted length functional $\int_0^l e^{\psi(\gamma(s))} \, |\gamma'(s)| \, ds$ in its relative homotopy class. For each $s \in [0,l]$, we denote by $\nu(s) \in T_{\gamma(s)} \Sigma$ the unit normal at the point $\gamma(s)$ and by $H(s)$ the geodesic curvature at the point $\gamma(s)$. The first variation formula implies $H(s) + \big \langle \nabla \psi|_{\gamma(s)},\nu(s) \big \rangle = 0$ for each $s \in [0,l]$. The stability inequality gives 
\begin{align*} 
&\int_0^l e^{\psi(\gamma(s))} \, \zeta'(s)^2 \, ds - \int_0^l e^{\psi(\gamma(s))} \, K(\gamma(s)) \, \zeta(s)^2 \, ds \\ 
&- \int_0^l e^{\psi(\gamma(s))} \, H(s)^2 \, \zeta(s)^2 \, ds + \int_0^l e^{\psi(\gamma(s))} \, (D^2 \psi)_{\gamma(s)}(\nu(s),\nu(s)) \, \zeta(s)^2 \, ds \\ 
&- e^{\psi(\gamma(0))} \, \kappa(\gamma(0)) \, \zeta(0)^2 - e^{\psi(\gamma(l))} \, \kappa(\gamma(l)) \, \zeta(l)^2 \geq 0 
\end{align*}
for every test function $\zeta \in C^\infty([0,l])$ (see Theorem \ref{second.variation}). We next consider the first eigenfunction of the stability operator. This gives a nonnegative function $v \in C^\infty([0,l])$ such that $\int_0^l e^{\psi(\gamma(s))} \, v(s)^2 \, ds = 1$ and 
\begin{align} 
\label{ode.for.v}
&-v''(s) - K(\gamma(s)) \, v(s) - H(s)^2 \, v(s) \notag \\ 
&+ (D^2 \psi)_{\gamma(s)}(\nu(s),\nu(s)) \, v(s) - \big \langle \nabla \psi|_{\gamma(s)},\gamma'(s) \big \rangle \, v'(s) = \lambda \, v(s) 
\end{align} 
for each $s \in [0,l]$, where $\lambda$ is a nonnegative constant. Moreover, $v$ satisfies the Neumann boundary conditions 
\begin{equation} 
\label{derivative.of.v.at.0}
-v'(0) = \kappa(\gamma(0)) \, v(0) 
\end{equation}
and 
\begin{equation} 
\label{derivative.of.v.at.l}
v'(l) = \kappa(\gamma(l)) \, v(l). 
\end{equation}
It is easy to see that $v(s) > 0$ for each $s \in [0,l]$. We next define 
\begin{equation} 
\label{definition.of.w}
w(s) = \psi(\gamma(s)) + \log v(s) 
\end{equation}
for each $s \in [0,l]$. 

\begin{lemma}
\label{ode.for.w} 
The function $w$ satisfies 
\[-w''(s) - \frac{N}{2(N-1)} \, w'(s)^2 + \frac{N(N-1)}{2} \geq 0\] 
for each $s \in [0,l]$. 
\end{lemma} 

\textbf{Proof.} 
Using (\ref{ode.for.v}), we obtain 
\begin{align}
\label{ode.for.log.v}
&-\frac{d^2}{ds^2} \log v(s) - v(s)^{-2} \, v'(s)^2 - K(\gamma(s)) - H(s)^2 \notag \\ 
&+ (D^2 \psi)_{\gamma(s)}(\nu(s),\nu(s)) - \big \langle \nabla \psi|_{\gamma(s)},\gamma'(s) \big \rangle \, v(s)^{-1} \, v'(s) = \lambda  
\end{align} 
for each $s \in [0,l]$. Moreover, 
\begin{align} 
\label{ode.for.psi}
&\Delta \psi(\gamma(s)) - (D^2 \psi)_{\gamma(s)}(\nu(s),\nu(s)) \notag \\ 
&- H(s) \, \big \langle \nabla \psi|_{\gamma(s)},\nu(s) \big \rangle - \frac{d^2}{ds^2} \psi(\gamma(s)) = 0
\end{align}
for each $s \in [0,l]$. In the next step, we add (\ref{ode.for.log.v}) and (\ref{ode.for.psi}). Using (\ref{definition.of.w}) and the identity $H(s) + \big \langle \nabla \psi|_{\gamma(s)},\nu(s) \big \rangle = 0$, we obtain 
\begin{align*} 
&-w''(s) + \Delta \psi(\gamma(s)) - K(\gamma(s)) \\ 
&- v(s)^{-2} \, v'(s)^2 - \big \langle \nabla \psi|_{\gamma(s)},\gamma'(s) \big \rangle \, v(s)^{-1} \, v'(s) = \lambda 
\end{align*} 
for each $s \in [0,l]$. Rearranging terms gives 
\begin{align*} 
&-w''(s) - \frac{N}{2(N-1)} \, w'(s)^2 \\ 
&+ \Delta \psi(\gamma(s)) + \frac{N-1}{2(N-2)} \, \big \langle \nabla \psi|_{\gamma(s)},\gamma'(s) \big \rangle^2 - K(\gamma(s)) \\ 
&- \frac{1}{2(N-2)(N-1)} \, \Big ( \big \langle \nabla \psi|_{\gamma(s)},\gamma'(s) \big \rangle - (N-2) \, v(s)^{-1} \, v'(s) \Big )^2 = \lambda 
\end{align*}
for each $s \in [0,l]$. By assumption, 
\[-\Delta \psi - \frac{N-1}{2(N-2)} \, |\nabla \psi|^2 + K + \frac{N(N-1)}{2} \geq 0\] 
at each point in $\Sigma$. Moreover, $\lambda$ is nonnegative. Putting these facts together, the assertion follows. This completes the proof of Lemma \ref{ode.for.w}. \\

\begin{lemma}
\label{derivative.of.w.at.0.and.l}
We have $-w'(0) > N-1$ and $w'(l) > N-1$.
\end{lemma}

\textbf{Proof.} 
Since the curve $\gamma$ meets $\partial \Sigma$ orthogonally, we know that $\eta|_{\gamma(0)} = -\gamma'(0)$ and $\eta|_{\gamma(l)} = \gamma'(l)$. Using (\ref{derivative.of.v.at.0}), (\ref{derivative.of.v.at.l}), and (\ref{definition.of.w}), we obtain 
\[-w'(0) = -\big \langle \nabla \psi|_{\gamma(0)},\gamma'(0) \big \rangle + \kappa(\gamma(0)) = \big \langle \nabla \psi|_{\gamma(0)},\eta|_{\gamma(0)} \big \rangle + \kappa(\gamma(0))\] 
and 
\[w'(l) = \big \langle \nabla \psi|_{\gamma(l)},\gamma'(l) \big \rangle + \kappa(\gamma(l)) = \big \langle \nabla \psi|_{\gamma(l)},\eta|_{\gamma(l)} \big \rangle + \kappa(\gamma(l)).\] 
On the other hand, it follows from (\ref{lower.bound.for.weighted.geodesic.curvature}) that $\langle \nabla \psi,\eta \rangle + \kappa > N-1$ at each point on $\partial \Sigma$. This completes the proof of Lemma \ref{derivative.of.w.at.0.and.l}. \\

Lemma \ref{derivative.of.w.at.0.and.l} implies that $-w'(0) > N-1$. Using Lemma \ref{ode.for.w} and standard ODE arguments, we conclude that $-w'(s) > N-1$ for each $s \in [0,l]$. In particular, $-w'(l) > N-1$. On the other hand, $w'(l) > N-1$ by Lemma \ref{derivative.of.w.at.0.and.l}. This is a contradiction. This completes the proof of Theorem \ref{2D.inequality} in the special case when $\partial \Sigma$ is disconnected.

\section{Proof of Theorem \ref{higher.dim.inequality.1}}

In this section, we explain how Theorem \ref{higher.dim.inequality.1} can be deduced from Theorem \ref{2D.inequality}. Let $n$ be an integer with $3 \leq n \leq 7$, let $N > n$ be a real number, let $M$ be a manifold of dimension $n$ with non-empty boundary $\partial M$, let $g$ be a Riemannian metric on $M$, and let $\varphi$ be a smooth function on $M$. We assume that $M$, $g$, and $\varphi$ satisfy the assumptions of Theorem \ref{higher.dim.inequality.1}. Let $\eta$ denote the outward-pointing unit normal vector field to $\partial M$. We denote by $h_{\partial M}$ the second fundamental form of $\partial M$ and by $H_{\partial M}$ the mean curvature of $\partial M$. Throughout this section, we assume that 
\begin{equation} 
\label{positivity.of.weighted.mean.curvature}
\inf_{\partial M} (\langle \nabla^M \varphi,\eta \rangle + H_{\partial M}) > 0, 
\end{equation}
for otherwise the assertion is trivial. 

Let $\Theta_0$ denote the pull-back of the volume form on $S^1$ under the map $\theta_0: \partial M \to S^1$. For each $k \in \{1,\hdots,n-2\}$, we denote by $\Theta_k$ the pull-back of the volume form on $S^1$ under the map $\theta_k|_{\partial M}: \partial M \to S^1$. Note that $\Theta_0,\Theta_1,\hdots,\Theta_{n-2}$ are closed one-forms on $\partial M$. By assumption, the map $(\theta_0,\theta_1|_{\partial M},\hdots,\theta_{n-2}|_{\partial M}): \partial M \to S^1 \times T^{n-2}$ has non-zero degree. This implies 
\[\int_{\partial M} \Theta_0 \wedge \Theta_1 \wedge \hdots \wedge \Theta_{n-2} \neq 0\] 
(see \cite{Lee}, Theorem 17.35 (a)).

\begin{proposition} 
\label{slicing}
We can find a collection of compact, connected, orientable submanifolds $\Sigma_k$, $k \in \{0,\hdots,n-2\}$, a collection of positive functions $\rho_k \in C^\infty(\Sigma_k)$, $k \in \{0,\hdots,n-2\}$, and a collection of positive functions $v_k \in C^\infty(\Sigma_k)$, $k \in \{1,\hdots,n-2\}$, with the following properties. 
\begin{itemize}
\item[(i)] $\Sigma_0 = M$ and $\rho_0 = e^\varphi$.
\item[(ii)] For each $k \in \{0,\hdots,n-2\}$, we have $\dim \Sigma_k = n-k$.
\item[(iii)] For each $k \in \{1,\hdots,n-2\}$, $\Sigma_k$ is a compact, connected, embedded, orientable hypersurface in $\Sigma_{k-1}$ satisfying $\partial \Sigma_k \subset \partial \Sigma_{k-1}$. Moreover, $\Sigma_k$ meets $\partial \Sigma_{k-1}$ orthogonally along $\partial \Sigma_k$. 
\item[(iv)] For each $k \in \{0,\hdots,n-2\}$, the outward-pointing unit normal vector field to $\partial \Sigma_k$ in $\Sigma_k$ equals $\eta$. Moreover, the second fundamental form of $\partial \Sigma_k$ in $\Sigma_k$ equals the restriction of $h_{\partial M}$ to $T(\partial \Sigma_k)$. 
\item[(v)] For each $k \in \{1,\hdots,n-2\}$, $\Sigma_k$ is a stable free boundary minimal hypersurface in $(\Sigma_{k-1},\rho_{k-1}^{\frac{2}{n-k}} \, g_{\Sigma_{k-1}})$.
\item[(vi)] For each $k \in \{1,\hdots,n-2\}$, the function $v_k \in C^\infty(\Sigma_k)$ satisfies 
\begin{align*} 
&-\Delta_{\Sigma_k} v_k - \text{\rm Ric}_{\Sigma_{k-1}}(\nu_{\Sigma_k},\nu_{\Sigma_k}) \, v_k - |h_{\Sigma_k}|^2 \, v_k \\ 
&+ (D_{\Sigma_{k-1}}^2 \log \rho_{k-1})(\nu_{\Sigma_k},\nu_{\Sigma_k}) \, v_k - \langle \nabla^{\Sigma_k} \log \rho_{k-1},\nabla^{\Sigma_k} v_k \rangle = \lambda_k v_k 
\end{align*}
on $\Sigma_k$ with Neumann boundary condition 
\[\langle \nabla^{\Sigma_k} v_k,\eta \rangle = h_{\partial M}(\nu_{\Sigma_k},\nu_{\Sigma_k}) \, v_k\] 
on $\partial \Sigma_k$. Here, $\lambda_k$ is a nonnegative constant.
\item[(vii)] For each $k \in \{1,\hdots,n-2\}$, the function $\rho_k \in C^\infty(\Sigma_k)$ is given by $\rho_k = \rho_{k-1}|_{\Sigma_k} \cdot v_k$.
\item[(viii)] For each $k \in \{0,\hdots,n-2\}$, the normal derivative of $\rho_k$ satisfies 
\[\langle \nabla^{\Sigma_k} \log \rho_k,\eta \rangle = \langle \nabla^M \varphi,\eta \rangle + \sum_{j=1}^k h_{\partial M}(\nu_{\Sigma_j},\nu_{\Sigma_j})\] 
at each point on $\partial \Sigma_k$.
\item[(ix)] For each $k \in \{0,\hdots,n-2\}$, the boundary mean curvature of the manifold $(\Sigma_k,\rho_k^{\frac{2}{n-k-1}} \, g_{\Sigma_k})$ is given by $\rho_k^{-\frac{1}{n-k-1}} \, (\langle \nabla^M \varphi,\eta \rangle + H_{\partial M})$ at each point on $\partial \Sigma_k$. 
\item[(x)] For each $k \in \{0,\hdots,n-2\}$, we have 
\[\int_{\partial \Sigma_k} \Theta_0 \wedge \Theta_1 \wedge \hdots \wedge \Theta_{n-k-2} \neq 0.\] 
\end{itemize}
\end{proposition}

\textbf{Proof.} 
We argue by induction on $k$. For $k=0$, we define $\Sigma_0 = M$ and $\rho_0 = e^\varphi$. It is clear that $\Sigma_0$ and $\rho_0$ satisfy all the required properties. We now turn to the inductive step. Suppose that $k \in \{1,\hdots,n-2\}$, and that we have constructed submanifolds $\Sigma_0,\hdots,\Sigma_{k-1}$, positive functions $\rho_0 \in C^\infty(\Sigma_0),\hdots,\rho_{k-1} \in C^\infty(\Sigma_{k-1})$, and positive functions $v_1 \in C^\infty(\Sigma_1),\hdots,v_{k-1} \in C^\infty(\Sigma_{k-1})$ satisfying the conditions (i)--(x) above. The inductive hypothesis implies that 
\[\int_{\partial \Sigma_{k-1}} \Theta_0 \wedge \Theta_1 \wedge \hdots \wedge \Theta_{n-k-1} \neq 0.\] 
By Sard's theorem, we can find an element $t \in S^1$ such that $t$ is a regular value of the function $\theta_{n-k-1}|_{\Sigma_{k-1}}: \Sigma_{k-1} \to S^1$ and $t$ is a regular value of the function $\theta_{n-k-1}|_{\partial \Sigma_{k-1}}: \partial \Sigma_{k-1} \to S^1$. We define $\tilde{\Sigma}_k = \Sigma_{k-1} \cap \{\theta_{n-k-1} = t\}$. Then $\tilde{\Sigma}_k$ is a compact, embedded, orientable hypersurface in $\Sigma_{k-1}$ with boundary $\partial \tilde{\Sigma}_k = \partial \Sigma_{k-1} \cap \{\theta_{n-k-1} = t\}$, and 
\[\int_{\partial \tilde{\Sigma}_k} \Theta_0 \wedge \Theta_1 \wedge \hdots \wedge \Theta_{n-k-2} \neq 0.\] 
Note that $\tilde{\Sigma}_k$ may be disconnected. This does not affect the subsequent arguments.

In view of the inductive hypothesis, the boundary mean curvature of the manifold $(\Sigma_{k-1},\rho_{k-1}^{\frac{2}{n-k}} \, g_{\Sigma_{k-1}})$ is given by $\rho_{k-1}^{-\frac{1}{n-k}} \, (\langle \nabla^M \varphi,\eta \rangle + H_{\partial M})$ at each point on $\partial \Sigma_{k-1}$. Using the inequality (\ref{positivity.of.weighted.mean.curvature}), we conclude that the manifold $(\Sigma_{k-1},\rho_{k-1}^{\frac{2}{n-k}} \, g_{\Sigma_{k-1}})$ has strictly mean convex boundary. Moreover, $\Theta_0 \wedge \Theta_1 \wedge \hdots \wedge \Theta_{n-k-2}$ is a closed $(n-k-1)$-form on $\partial \Sigma_{k-1}$. By Theorem \ref{existence.and.regularity}, we can find a compact, connected, embedded, orientable hypersurface $\Sigma_k$ in $\Sigma_{k-1}$ with the following properties: 
\begin{itemize}
\item The boundary $\partial \Sigma_k$ is contained in $\partial \Sigma_{k-1}$. Moreover, $\Sigma_k$ meets $\partial \Sigma_{k-1}$ orthogonally along $\partial \Sigma_k$. 
\item The submanifold $\Sigma_k$ is a stable free boundary minimal hypersurface in $(\Sigma_{k-1},\rho_{k-1}^{\frac{2}{n-k}} \, g_{\Sigma_{k-1}})$. 
\item We have $\int_{\partial \Sigma_k} \Theta_0 \wedge \Theta_1 \wedge \hdots \wedge \Theta_{n-k-2} \neq 0$. 
\end{itemize}
In view of the inductive hypothesis, the outward-pointing unit normal vector field to $\partial \Sigma_{k-1}$ in $\Sigma_{k-1}$ equals $\eta$. Consequently, the outward-pointing unit normal vector field to $\partial \Sigma_k$ in $\Sigma_k$ equals $\eta$. From this, we deduce that the second fundamental form of $\partial \Sigma_k$ in $\Sigma_k$ equals the restriction of $h_{\partial M}$ to $T(\partial \Sigma_k)$. In particular, the mean curvature of $\partial \Sigma_k$ in $\Sigma_k$ is given by 
\begin{equation} 
\label{boundary.mean.curvature.of.Sigma_k}
\text{\rm tr}_{\partial \Sigma_k}(h_{\partial M}) = H_{\partial M} - \sum_{j=1}^k h_{\partial M}(\nu_{\Sigma_j},\nu_{\Sigma_j}). 
\end{equation}
To summarize, we have shown that properties (ii), (iii), (iv), (v), and (x) hold for $\Sigma_k$. 

The stability inequality implies that 
\begin{align*} 
&\int_{\Sigma_k} \rho_{k-1} \, |\nabla^{\Sigma_k} \zeta|^2 - \int_{\Sigma_k} \rho_{k-1} \, \text{\rm Ric}_{\Sigma_{k-1}}(\nu_{\Sigma_k},\nu_{\Sigma_k}) \, \zeta^2 \\ 
&- \int_{\Sigma_k} \rho_{k-1} \, |h_{\Sigma_k}|^2 \, \zeta^2 + \int_{\Sigma_k} \rho_{k-1} \, (D_{\Sigma_{k-1}}^2 \log \rho_{k-1})(\nu_{\Sigma_k},\nu_{\Sigma_k}) \, \zeta^2 \\ 
&- \int_{\partial \Sigma_k} \rho_{k-1} \, h_{\partial \Sigma_{k-1}}(\nu_{\Sigma_k},\nu_{\Sigma_k}) \, \zeta^2 \geq 0 
\end{align*}
for every test function $\zeta \in C^\infty(\Sigma_k)$ (see Theorem \ref{second.variation} below). Here, $h_{\partial \Sigma_{k-1}}$ denotes the second fundamental form of $\partial \Sigma_{k-1}$ in $\Sigma_{k-1}$. The inductive hypothesis implies that $h_{\partial \Sigma_{k-1}}(\nu_{\Sigma_k},\nu_{\Sigma_k}) =  h_{\partial M}(\nu_{\Sigma_k},\nu_{\Sigma_k})$. We next consider the first eigenfunction of the stability operator. This gives a nonnegative function $v_k \in C^\infty(\Sigma_k)$ such that $\int_{\Sigma_k} \rho_{k-1} \, v_k^2 = 1$ and 
\begin{align*} 
&-\Delta_{\Sigma_k} v_k - \text{\rm Ric}_{\Sigma_{k-1}}(\nu_{\Sigma_k},\nu_{\Sigma_k}) \, v_k - |h_{\Sigma_k}|^2 \, v_k \\ 
&+ (D_{\Sigma_{k-1}}^2 \log \rho_{k-1})(\nu_{\Sigma_k},\nu_{\Sigma_k}) \, v_k - \langle \nabla^{\Sigma_k} \log \rho_{k-1},\nabla^{\Sigma_k} v_k \rangle = \lambda_k v_k 
\end{align*}
at each point on $\Sigma_k$, where $\lambda_k$ is a nonnegative constant. Moreover, $v_k$ satisfies the Neumann boundary condition 
\begin{equation} 
\label{normal.derivative.of.v_k}
\langle \nabla^{\Sigma_k} v_k,\eta \rangle = h_{\partial M}(\nu_{\Sigma_k},\nu_{\Sigma_k}) \, v_k 
\end{equation}
at each point on $\partial \Sigma_k$. Using the strict maximum principle and the Hopf boundary point lemma, we conclude that $v_k$ is strictly positive at each point on $\Sigma_k$. Therefore, property (vi) holds for $v_k$. 

We next define the function $\rho_k \in C^\infty(\Sigma_k)$ by $\rho_k = \rho_{k-1}|_{\Sigma_k} \cdot v_k$. Then property (vii) holds for $\rho_k$. The inductive hypothesis implies that 
\begin{equation} 
\label{normal.derivative.of.log.rho_k-1}
\langle \nabla^{\Sigma_{k-1}} \log \rho_{k-1},\eta \rangle = \langle \nabla^M \varphi,\eta \rangle + \sum_{j=1}^{k-1} h_{\partial M}(\nu_{\Sigma_j},\nu_{\Sigma_j}) 
\end{equation}
at each point on $\partial \Sigma_{k-1}$. Combining (\ref{normal.derivative.of.v_k}) and (\ref{normal.derivative.of.log.rho_k-1}), we conclude that 
\begin{equation} 
\label{normal.derivative.of.log.rho_k}
\langle \nabla^{\Sigma_k} \log \rho_k,\eta \rangle = \langle \nabla^M \varphi,\eta \rangle + \sum_{j=1}^k h_{\partial M}(\nu_{\Sigma_j},\nu_{\Sigma_j}) 
\end{equation}
at each point on $\partial \Sigma_k$. Therefore, property (viii) holds for $\rho_k$.

Finally, combining (\ref{boundary.mean.curvature.of.Sigma_k}) and (\ref{normal.derivative.of.log.rho_k}), we conclude that the boundary mean curvature of the manifold $(\Sigma_k,\rho_k^{\frac{2}{n-k-1}} \, g_{\Sigma_k})$ is given by 
\begin{align*} 
&\rho_k^{-\frac{1}{n-k-1}} \, \Big ( \langle \nabla^{\Sigma_k} \log \rho_k,\eta \rangle + H_{\partial M} - \sum_{j=1}^k h_{\partial M}(\nu_{\Sigma_j},\nu_{\Sigma_j}) \Big ) \\ 
&= \rho_k^{-\frac{1}{n-k-1}} \, (\langle \nabla^M \varphi,\eta \rangle + H_{\partial M}) 
\end{align*} 
at each point on $\partial \Sigma_k$. Thus, property (ix) holds for $\Sigma_k$. This completes the proof of Proposition \ref{slicing}. \\

We next state an identity that relates the scalar curvature of $\Sigma_k$ to the scalar curvature of $\Sigma_{k-1}$. This identity is similar to a formula used by Schoen and Yau \cite{Schoen-Yau1},\cite{Schoen-Yau3} in their dimension descent scheme. It is closely related to the torical symmetrization procedure described in Theorem 11.1 of Gromov and Lawson's paper \cite{Gromov-Lawson} (see also \cite{Fischer-Colbrie-Schoen} and \cite{Gromov}).

\begin{proposition}[cf. R.~Schoen, S.T.~Yau \cite{Schoen-Yau1},\cite{Schoen-Yau3}; M.~Gromov, H.B.~Lawson, Jr. \cite{Gromov-Lawson}]
\label{scalar.curvature.of.consecutive.slices}
Let $k \in \{1,\hdots,n-2\}$. Then 
\begin{align*} 
&- 2 \, \Delta_{\Sigma_k} \log \rho_k - |\nabla^{\Sigma_k} \log \rho_k|^2 + R_{\Sigma_k} \\ 
&+ 2 \, \Delta_{\Sigma_{k-1}} \log \rho_{k-1} + |\nabla^{\Sigma_{k-1}} \log \rho_{k-1}|^2 - R_{\Sigma_{k-1}} - |\nabla^{\Sigma_k} \log v_k|^2 - |h_{\Sigma_k}|^2 \\ 
&= 2 \lambda_k 
\end{align*} 
at each point on $\Sigma_k$.
\end{proposition}

\textbf{Proof.} 
Using the Gauss equations, we obtain 
\begin{equation} 
\label{gauss.equation}
R_{\Sigma_k} - R_{\Sigma_{k-1}} + 2 \,  \text{\rm Ric}_{\Sigma_{k-1}}(\nu_{\Sigma_k},\nu_{\Sigma_k}) - H_{\Sigma_k}^2 + |h_{\Sigma_k}|^2 = 0 
\end{equation}
at each point on $\Sigma_k$. Property (vi) in Proposition \ref{slicing} implies 
\begin{align} 
\label{pde.for.log.v_k}
&-2 \, \Delta_{\Sigma_k} \log v_k - 2 \, |\nabla^{\Sigma_k} \log v_k|^2 - 2 \, \text{\rm Ric}_{\Sigma_{k-1}}(\nu_{\Sigma_k},\nu_{\Sigma_k}) - 2 \, |h_{\Sigma_k}|^2 \notag \\ 
&+ 2 \, (D_{\Sigma_{k-1}}^2 \log \rho_{k-1})(\nu_{\Sigma_k},\nu_{\Sigma_k}) - 2 \, \langle \nabla^{\Sigma_k} \log \rho_{k-1},\nabla^{\Sigma_k} \log v_k \rangle = 2\lambda_k 
\end{align} 
at each point on $\Sigma_k$. Moreover, 
\begin{align} 
\label{laplacian}
&2 \, \Delta_{\Sigma_{k-1}} \log \rho_{k-1} - 2 \, (D_{\Sigma_{k-1}}^2 \log \rho_{k-1})(\nu_{\Sigma_k},\nu_{\Sigma_k}) \notag \\ 
&- 2 \, H_{\Sigma_k} \, \langle \nabla^{\Sigma_{k-1}} \log \rho_{k-1},\nu_{\Sigma_k} \rangle - 2 \, \Delta_{\Sigma_k} \log \rho_{k-1} = 0 
\end{align} 
at each point on $\Sigma_k$. By property (v) in Proposition \ref{slicing}, $\Sigma_k$ is a minimal hypersurface in $(\Sigma_{k-1},\rho_{k-1}^{\frac{2}{n-k}} \, g_{\Sigma_{k-1}})$. This implies 
\begin{equation} 
\label{first.variation}
H_{\Sigma_k} + \langle \nabla^{\Sigma_{k-1}} \log \rho_{k-1},\nu_{\Sigma_k} \rangle = 0 
\end{equation}
at each point on $\Sigma_k$. Using (\ref{first.variation}) and the identity 
\[|\nabla^{\Sigma_{k-1}} \log \rho_{k-1}|^2 - |\nabla^{\Sigma_k} \log \rho_{k-1}|^2 = \langle \nabla^{\Sigma_{k-1}} \log \rho_{k-1},\nu_{\Sigma_k} \rangle^2,\] 
we obtain 
\begin{align} 
\label{consequence.of.first.variation}
&H_{\Sigma_k}^2 + 2 \, H_{\Sigma_k} \, \langle \nabla^{\Sigma_{k-1}} \log \rho_{k-1},\nu_{\Sigma_k} \rangle \notag \\ 
&+ |\nabla^{\Sigma_{k-1}} \log \rho_{k-1}|^2 - |\nabla^{\Sigma_k} \log \rho_{k-1}|^2 = 0 
\end{align}
at each point on $\Sigma_k$. In the next step, we add (\ref{gauss.equation}), (\ref{pde.for.log.v_k}), (\ref{laplacian}), and (\ref{consequence.of.first.variation}). This gives 
\begin{align*} 
&- 2 \, (\Delta_{\Sigma_k} \log \rho_{k-1} + \Delta_{\Sigma_k} \log v_k) - |\nabla^{\Sigma_k} \log \rho_{k-1} + \nabla^{\Sigma_k} \log v_k|^2 + R_{\Sigma_k} \\ 
&+ 2 \, \Delta_{\Sigma_{k-1}} \log \rho_{k-1} + |\nabla^{\Sigma_{k-1}} \log \rho_{k-1}|^2 - R_{\Sigma_{k-1}} - |\nabla^{\Sigma_k} \log v_k|^2 - |h_{\Sigma_k}|^2 \\ 
&= 2 \lambda_k 
\end{align*}
at each point on $\Sigma_k$. Finally, $\log \rho_k = \log \rho_{k-1} + \log v_k$ at each point on $\Sigma_k$. Putting these facts together, the assertion follows. This completes the proof of Proposition \ref{scalar.curvature.of.consecutive.slices}. \\

\begin{corollary}
\label{scalar.curvature.of.consecutive.slices.2}
Let $k \in \{1,\hdots,n-2\}$. Then 
\begin{align*} 
&- 2 \, \Delta_{\Sigma_k} \log \rho_k - \frac{N-n+k+1}{N-n+k} \, |\nabla^{\Sigma_k} \log \rho_k|^2 + R_{\Sigma_k} \\ 
&+ 2 \, \Delta_{\Sigma_{k-1}} \log \rho_{k-1} + \frac{N-n+k}{N-n+k-1} \, |\nabla^{\Sigma_{k-1}} \log \rho_{k-1}|^2 - R_{\Sigma_{k-1}} \geq 0 
\end{align*} 
at each point on $\Sigma_k$.
\end{corollary}

\textbf{Proof.} 
Using Proposition \ref{scalar.curvature.of.consecutive.slices} and the inequality $\lambda_k \geq 0$, we obtain 
\begin{align*} 
&- 2 \, \Delta_{\Sigma_k} \log \rho_k - |\nabla^{\Sigma_k} \log \rho_k|^2 + R_{\Sigma_k} \\ 
&+ 2 \, \Delta_{\Sigma_{k-1}} \log \rho_{k-1} + |\nabla^{\Sigma_{k-1}} \log \rho_{k-1}|^2 - R_{\Sigma_{k-1}} - |\nabla^{\Sigma_k} \log v_k|^2 \geq 0 
\end{align*} 
at each point on $\Sigma_k$. It follows from property (vii) in Proposition \ref{slicing} that $\rho_k = \rho_{k-1} \, v_k$ at each point on $\Sigma_k$. This implies 
\begin{align*} 
&\frac{1}{N-n+k-1} \, |\nabla^{\Sigma_{k-1}} \log \rho_{k-1}|^2 - \frac{1}{N-n+k} \, |\nabla^{\Sigma_k} \log \rho_k|^2 + |\nabla^{\Sigma_k} \log v_k|^2 \\ 
&= \frac{N-n+k-1}{N-n+k} \, \Big | \frac{1}{N-n+k-1} \, \nabla^{\Sigma_k} \log \rho_{k-1} - \nabla^{\Sigma_k} \log v_k \Big |^2 \\ 
&+ \frac{1}{N-n+k-1} \, \langle \nabla^{\Sigma_{k-1}} \log \rho_{k-1},\nu_{\Sigma_k} \rangle^2 \geq 0
\end{align*} 
at each point on $\Sigma_k$. Adding these two inequalities, the assertion follows. This completes the proof of Corollary \ref{scalar.curvature.of.consecutive.slices.2}. \\

\begin{corollary}
\label{inequality.for.scalar.curvature}
Let $k \in \{0,\hdots,n-2\}$. Then 
\begin{align*} 
&-2 \, \Delta_{\Sigma_k} \log \rho_k - \frac{N-n+k+1}{N-n+k} \, |\nabla^{\Sigma_k} \log \rho_k|^2 + R_{\Sigma_k} \\ 
&+ 2 \, \Delta_M \varphi + \frac{N-n+1}{N-n} \, |\nabla^M \varphi|^2 - R_M \geq 0 
\end{align*}
at each point on $\Sigma_k$.
\end{corollary}

\textbf{Proof.} 
The proof is by induction on $k$. For $k=0$, the assertion is trivial. The inductive step follows from Corollary \ref{scalar.curvature.of.consecutive.slices.2}. This completes the proof of Corollary \ref{inequality.for.scalar.curvature}. \\

After these preparations, we now complete the proof of Theorem \ref{higher.dim.inequality.1}. To that end, we consider the two-dimensional surface $\Sigma = \Sigma_{n-2}$. We define a function $\psi \in C^\infty(\Sigma)$ by $\psi = \log \rho_{n-2}$. It follows from property (x) in Proposition \ref{slicing} that $\int_{\partial \Sigma} \Theta_0 \neq 0$. In particular, $\Sigma$ has non-empty boundary. Let $K$ denote the Gaussian curvature of $\Sigma$, and let $\kappa$ denote the geodesic curvature of the boundary $\partial \Sigma$. 

By assumption, 
\[-2 \, \Delta_M \varphi - \frac{N-n+1}{N-n} \, |\nabla^M \varphi|^2 + R_M + N(N-1) \geq 0\] 
at each point in $M$. On the other hand, applying Corollary \ref{inequality.for.scalar.curvature} with $k=n-2$ gives 
\begin{align*} 
&-2 \, \Delta \psi - \frac{N-1}{N-2} \, |\nabla \psi|^2 + 2K \\ 
&+ 2 \, \Delta_M \varphi + \frac{N-n+1}{N-n} \, |\nabla^M \varphi|^2 - R_M \geq 0 
\end{align*}
at each point on $\Sigma$. Adding these inequalities, we obtain 
\[-2 \, \Delta \psi - \frac{N-1}{N-2} \, |\nabla \psi|^2 + 2K + N(N-1) \geq 0\] 
at each point on $\Sigma$. Using Theorem \ref{2D.inequality}, we conclude that 
\[2 \, |\partial \Sigma|^N \, \inf_{\partial \Sigma} (\langle \nabla \psi,\eta \rangle + \kappa - (N-1)) < \Big ( \frac{4\pi}{N} \Big )^N.\] 
By property (ix) in Proposition \ref{slicing}, the geodesic curvature of $\partial \Sigma$ with respect to the conformal metric $e^{2\psi} \, g_\Sigma$ is given by $e^{-\psi} \, (\langle \nabla^M \varphi,\eta \rangle + H_{\partial M})$. In other words, $\langle \nabla \psi,\eta \rangle + \kappa = \langle \nabla^M \varphi,\eta \rangle + H_{\partial M}$ at each point on $\partial \Sigma$. Thus, 
\[2 \, |\partial \Sigma|^N \, \inf_{\partial \Sigma} (\langle \nabla^M \varphi,\eta \rangle + H_{\partial M} - (N-1)) < \Big ( \frac{4\pi}{N} \Big )^N.\] 
Finally, since $\int_{\partial \Sigma} \Theta_0 \neq 0$, it follows that $|\partial \Sigma| \geq \sigma$ by definition of $\sigma$. This completes the proof of Theorem \ref{higher.dim.inequality.1}.

\section{Proof of Theorem \ref{higher.dim.inequality.2}}

In this section, we give the proof of Theorem \ref{higher.dim.inequality.2}. Let $n$ be an integer with $3 \leq n \leq 7$, let $M$ be a manifold of dimension $n$ with non-empty boundary $\partial M$, and let $g$ be a Riemannian metric on $M$. We assume that $M$ and $g$ satisfy the assumptions of Theorem \ref{higher.dim.inequality.2}. Throughout this section, we assume that 
\[\inf_{\partial M} H_{\partial M} > 0,\] 
for otherwise the assertion is trivial. 

We apply Proposition \ref{slicing} with $\varphi=0$. Hence, we can find a collection of compact, connected, orientable submanifolds $\Sigma_k$, $k \in \{0,\hdots,n-2\}$, a collection of positive functions $\rho_k \in C^\infty(\Sigma_k)$, $k \in \{0,\hdots,n-2\}$, and a collection of positive functions $v_k \in C^\infty(\Sigma_k)$, $k \in \{1,\hdots,n-2\}$, satisfying the conditions (i)--(x) in Proposition \ref{slicing}. We define $\Sigma = \Sigma_{n-2}$ and $\psi = \log \rho_{n-2}$. It follows from property (x) in Proposition \ref{slicing} that $\int_{\partial \Sigma} \Theta_0 \neq 0$. In particular, $\Sigma$ has non-empty boundary. 

By assumption, $R_M + n(n-1) \geq 0$ at each point in $M$. In the next step, we apply Corollary \ref{inequality.for.scalar.curvature} with $\varphi=0$ and $N > n$, and pass to the limit as $N \searrow n$. Hence, if $k \in \{1,\hdots,n-2\}$, then 
\[-2 \, \Delta_{\Sigma_k} \log \rho_k - \frac{k+1}{k} \, |\nabla^{\Sigma_k} \log \rho_k|^2 + R_{\Sigma_k} + n(n-1) \geq 0\] 
at each point on $\Sigma_k$. Putting $k = n-2$, we obtain 
\[-2 \, \Delta \psi - \frac{n-1}{n-2} \, |\nabla \psi|^2 + 2K + n(n-1) \geq 0\] 
at each point on $\Sigma$. Using Theorem \ref{2D.inequality}, we conclude that 
\[2 \, |\partial \Sigma|^n \, \inf_{\partial \Sigma} (\langle \nabla \psi,\eta \rangle + \kappa - (n-1)) < \Big ( \frac{4\pi}{n} \Big )^n.\] 
We now apply property (ix) in Proposition \ref{slicing} with $\varphi=0$. This implies $\langle \nabla \psi,\eta \rangle + \kappa = H_{\partial M}$ at each point on $\partial \Sigma$. Thus, 
\[2 \, |\partial \Sigma|^n \, \inf_{\partial \Sigma} (H_{\partial M} - (n-1)) < \Big ( \frac{4\pi}{n} \Big )^n.\] 
Finally, since $\int_{\partial \Sigma} \Theta_0 \neq 0$, it follows that $|\partial \Sigma| \geq \sigma$ by definition of $\sigma$. This completes the proof of Theorem \ref{higher.dim.inequality.2}.

\section{Proof of Theorem \ref{higher.dim.inequality.3}}

In this section, we give the proof of Theorem \ref{higher.dim.inequality.3}. Let $n$ be an integer with $3 \leq n \leq 7$, let $M$ be a manifold of dimension $n$ with non-empty boundary $\partial M$, let $g$ be a Riemannian metric on $M$, and let $\varphi$ be a smooth function on $M$. We assume that $M$, $g$, and $\varphi$ satisfy the assumptions of Theorem \ref{higher.dim.inequality.3}. Throughout this section, we assume that 
\[\inf_{\partial M} (\langle \nabla^M \varphi,\eta \rangle + H_{\partial M}) > 0,\] 
for otherwise the assertion is trivial. 

We can find a collection of compact, connected, orientable submanifolds $\Sigma_k$, $k \in \{0,\hdots,n-2\}$, a collection of positive functions $\rho_k \in C^\infty(\Sigma_k)$, $k \in \{0,\hdots,n-2\}$, and a collection of positive functions $v_k \in C^\infty(\Sigma_k)$, $k \in \{1,\hdots,n-2\}$, satisfying the conditions (i)--(x) in Proposition \ref{slicing}. We define $\Sigma = \Sigma_{n-2}$ and $\psi = \log \rho_{n-2}$. Property (x) in Proposition \ref{slicing} implies that $\int_{\partial \Sigma} \Theta_0 \neq 0$. In particular, $\Sigma$ has non-empty boundary. Since $\Sigma$ is connected, it follows that $\chi(\Sigma) \leq 1$. 

By assumption, 
\[-2 \, \Delta_M \varphi - |\nabla^M \varphi|^2 + R_M \geq 0\] 
at each point on $\Sigma$. For each $k \in \{1,\hdots,n-2\}$, Proposition \ref{scalar.curvature.of.consecutive.slices} implies that 
\begin{align*} 
&- 2 \, \Delta_{\Sigma_k} \log \rho_k - |\nabla^{\Sigma_k} \log \rho_k|^2 + R_{\Sigma_k} \\ 
&+ 2 \, \Delta_{\Sigma_{k-1}} \log \rho_{k-1} + |\nabla^{\Sigma_{k-1}} \log \rho_{k-1}|^2 - R_{\Sigma_{k-1}} \geq 0 
\end{align*} 
at each point on $\Sigma_k$. Hence, if $k \in \{1,\hdots,n-2\}$, then 
\[- 2 \, \Delta_{\Sigma_k} \log \rho_k - |\nabla^{\Sigma_k} \log \rho_k|^2 + R_{\Sigma_k} \geq 0\] 
at each point on $\Sigma_k$. Putting $k = n-2$, we obtain 
\begin{equation} 
\label{pde.for.psi}
-2 \, \Delta \psi - |\nabla \psi|^2 + 2K \geq 0 
\end{equation}
at each point on $\Sigma$. Since $\chi(\Sigma) \leq 1$, the Gauss-Bonnet theorem gives 
\begin{equation} 
\label{gauss.bonnet}
\int_{\partial \Sigma} \kappa \leq 2\pi - \int_\Sigma K. 
\end{equation}
Moreover, 
\begin{equation} 
\label{divergence.theorem}
\int_{\partial \Sigma} \langle \nabla \psi,\eta \rangle = \int_\Sigma \Delta \psi 
\end{equation} 
by the divergence theorem. In the next step, we add (\ref{gauss.bonnet}) and (\ref{divergence.theorem}). Using (\ref{pde.for.psi}), we conclude that 
\[\int_{\partial \Sigma} (\langle \nabla \psi,\eta \rangle + \kappa) \leq 2\pi + \int_\Sigma (\Delta \psi - K) \leq 2\pi.\] 
Consequently, 
\[|\partial \Sigma| \, \inf_{\partial \Sigma} (\langle \nabla \psi,\eta \rangle + \kappa) \leq 2\pi.\] 
Property (ix) in Proposition \ref{slicing} implies that $\langle \nabla \psi,\eta \rangle + \kappa = \langle \nabla^M \varphi,\eta \rangle + H_{\partial M}$ at each point on $\partial \Sigma$. Thus, 
\[|\partial \Sigma| \, \inf_{\partial \Sigma} (\langle \nabla^M \varphi,\eta \rangle + H_{\partial M}) \leq 2\pi.\] 
Finally, since $\int_{\partial \Sigma} \Theta_0 \neq 0$, it follows that $|\partial \Sigma| \geq \sigma$ by definition of $\sigma$. This completes the proof of Theorem \ref{higher.dim.inequality.3}.

\section{Proof of Theorem \ref{positive.energy.theorem}} 

In this section, we explain how Theorem \ref{positive.energy.theorem} follows from Theorem \ref{higher.dim.inequality.2}. Let $n$ be an integer with $3 \leq n \leq 7$, let $\gamma$ be a flat metric on $S^1 \times T^{n-2}$, let $Q$ be a symmetric $(0,2)$-tensor on $S^1 \times T^{n-2}$, let $M$ be a manifold of dimension $n$, and let $g$ be a Riemannian metric on $M$. Throughout this section, we assume that the assumptions of Theorem \ref{positive.energy.theorem} are satisfied. We define a real number $\mu$ by 
\begin{equation} 
\label{definition.of.mu}
\int_{S^1 \times T^{n-2}} \big ( n \, \text{\rm tr}_\gamma(Q) + 2\mu \big ) \, d\text{\rm vol}_\gamma = 0. 
\end{equation} 
Let us fix a smooth function $u: S^1 \times T^{n-2} \to \mathbb{R}$ such that 
\begin{equation} 
\label{pde.for.u}
\Delta_\gamma u + \frac{n}{2} \, \text{\rm tr}_\gamma(Q) + \mu = 0 
\end{equation}
at each point on $S^1 \times T^{n-2}$. The function $u$ is unique up to additive constants. If we normalize $u$ so that $\int_{S^1 \times T^{n-2}} u \, d\text{\rm vol}_\gamma = 0$, then $u$ is uniquely determined. We denote by $\hat{u}: [r_0,\infty) \times S^1 \times T^{n-2} \to \mathbb{R}$ the composition of $u$ with the canonical projection from $[r_0,\infty) \times S^1 \times T^{n-2}$ to $S^1 \times T^{n-2}$. 

We define a hyperbolic metric $\bar{g}$ on $[r_0,\infty) \times S^1 \times T^{n-2}$ by 
\[\bar{g} = r^{-2} \, dr \otimes dr + r^2 \, \gamma.\] 
Moreover, we define a metric $\hat{g}$ on $[r_0,\infty) \times S^1 \times T^{n-2}$ by 
\[\hat{g} = \bar{g} + r^{2-n} \, Q = r^{-2} \, dr \otimes dr + r^2 \, \gamma + r^{2-n} \, Q.\] 
By assumption, $|g-\hat{g}|_{\bar{g}} \leq o(r^{-n})$ and $|\bar{D}(g-\hat{g})|_{\bar{g}} \leq o(r^{-n})$. 

\begin{lemma}
\label{Hessian.of.hat.u}
Let $D_\gamma^2 u$ denote the Hessian of the function $u: S^1 \times T^{n-2} \to \mathbb{R}$ with respect to the flat metric $\gamma$. Moreover, let $D^2 \hat{u}$ denote the Hessian of the function $\hat{u}: [r_0,\infty) \times S^1 \times T^{n-2} \to \mathbb{R}$ with respect to the metric $g$. Then 
\[\Big | D^2 \hat{u} - D_\gamma^2 u + r^{-1} \, (dr \otimes d\hat{u} + d\hat{u} \otimes dr) \Big |_g \leq O(r^{-n-1}).\] 
\end{lemma}

\textbf{Proof.} 
Let $\bar{D}^2 \hat{u}$ denote the Hessian of the function $\hat{u}$ with respect to the hyperbolic metric $\bar{g}$. A straightforward calculation gives 
\[\bar{D}^2 \hat{u} - D_\gamma^2 u + r^{-1} \, (dr \otimes d\hat{u} + d\hat{u} \otimes dr) = 0.\] 
Using the estimates $|g-\bar{g}|_{\bar{g}} \leq O(r^{-n})$ and $|\bar{D}(g-\bar{g})|_{\bar{g}} \leq O(r^{-n})$, we obtain 
\[|D^2 \hat{u} - \bar{D}^2 \hat{u}|_{\bar{g}} \leq C \, |\bar{D}(g-\bar{g})|_{\bar{g}} \, |d\hat{u}|_{\bar{g}} \leq O(r^{-n-1}).\] 
Putting these facts together, the assertion follows. This completes the proof of Lemma \ref{Hessian.of.hat.u}. \\

\begin{lemma} 
\label{Hessian.of.r}
Let $D^2 r$ denote the Hessian of the function $r$ with respect to the metric $g$. Then 
\[\Big | D^2 r - r \, g + \frac{n}{2} \, r^{3-n} \, Q \Big |_g \leq o(r^{1-n}).\] 
\end{lemma} 

\textbf{Proof.} 
Let $V = r^2 \, \frac{\partial}{\partial r}$. Then $\hat{g}(V,\cdot) = dr$. Therefore, the gradient of the function $r$ with respect to the metric $\hat{g}$ is equal to $V$. The Hessian of the function $r$ with respect to the metric $\hat{g}$ is given by $\hat{D}^2 r = \frac{1}{2} \, \mathscr{L}_V(\hat{g})$. A straightforward calculation shows that 
\[\mathscr{L}_V(\hat{g}) - 2r \, \hat{g} + n \, r^{3-n} \, Q = 0.\] 
This implies 
\[\hat{D}^2 r - r \, \hat{g} + \frac{n}{2} \, r^{3-n} \, Q = 0.\] 
Using the estimates $|g-\hat{g}|_{\bar{g}} \leq o(r^{-n})$ and $|\bar{D}(g-\hat{g})|_{\bar{g}} \leq o(r^{-n})$, we obtain $|\hat{D}(g-\hat{g})|_{\hat{g}} \leq o(r^{-n})$. Thus, we conclude that 
\[|D^2 r - \hat{D}^2 r|_{\hat{g}} \leq C \, |\hat{D}(g-\hat{g})|_{\hat{g}} \, |dr|_{\hat{g}} \leq o(r^{1-n}).\] 
Putting these facts together, the assertion follows. This completes the proof of Lemma \ref{Hessian.of.r}. \\

In the following, we consider a sequence $r_j \to \infty$. For $j$ sufficiently large, we define $M^{(j)} = M \setminus \{r > r_j + r_j^{3-n} \, \hat{u}\}$. Note that $M^{(j)}$ is a compact domain in $M$ with boundary $\partial M^{(j)} = \{r = r_j + r_j^{3-n} \, \hat{u}\}$.

\begin{proposition}
\label{asymptotic.formula.for.mean.curvature}
The mean curvature of the boundary $\partial M^{(j)}$ with respect to the metric $g$ satisfies 
\[\sup_{\partial M^{(j)}} \big | H_{\partial M^{(j)}} - (n-1) - r^{-n} \, \mu \big | \leq o(r_j^{-n}).\] 
\end{proposition}

\textbf{Proof.} 
The boundary $\partial M^{(j)}$ is a level set of the function $r - r_j^{3-n} \, \hat{u}$. The standard formula for the mean curvature of a level set gives 
\begin{equation} 
\label{mean.curvature} 
H_{\partial M^{(j)}} = \frac{\text{\rm tr}_{\partial M^{(j)}}(D^2 r) - r_j^{3-n} \, \text{\rm tr}_{\partial M^{(j)}}(D^2 \hat{u})}{|dr - r_j^{3-n} \, d\hat{u}|_g} 
\end{equation} 
at each point on $\partial M^{(j)}$. Here, $D^2 r$ denotes the Hessian of the function $r$ with respect to the metric $g$ and $D^2 \hat{u}$ denotes the Hessian of the function $\hat{u}$ with respect to the metric $g$. The boundary trace $\text{\rm tr}_{\partial M^{(j)}}$ is computed using the metric $g$. 

Using Lemma \ref{Hessian.of.r}, we obtain 
\[\Big | D^2 r - r \, g + \frac{n}{2} \, r^{3-n} \, Q \Big |_g \leq o(r_j^{1-n})\] 
at each point on $\partial M^{(j)}$. Using the identity $\text{\rm tr}_{\partial M^{(j)}}(g) = n-1$, it follows that 
\[\Big | \text{\rm tr}_{\partial M^{(j)}}(D^2 r) - (n-1) \, r + \frac{n}{2} \, r^{3-n} \, \text{\rm tr}_{\partial M^{(j)}}(Q) \Big | \leq o(r_j^{1-n})\] 
at each point on $\partial M^{(j)}$. Consequently, 
\begin{equation} 
\label{first.term}
\Big | \text{\rm tr}_{\partial M^{(j)}}(D^2 r) - (n-1) \, r + \frac{n}{2} \, r^{1-n} \, \text{\rm tr}_\gamma(Q) \Big | \leq o(r_j^{1-n}) 
\end{equation}
at each point on $\partial M^{(j)}$. Using Lemma \ref{Hessian.of.hat.u}, we obtain 
\[\Big | D^2 \hat{u} - D_\gamma^2 u + r^{-1} \, \big ( (dr - r_j^{3-n} \, d\hat{u}) \otimes d\hat{u} + d\hat{u} \otimes (dr - r_j^{3-n} \, d\hat{u}) \big ) \Big |_g \leq o(r_j^{-2})\] 
at each point on $\partial M^{(j)}$. Using the identity 
\[\text{\rm tr}_{\partial M^{(j)}} \big ( (dr - r_j^{3-n} \, d\hat{u}) \otimes d\hat{u} + d\hat{u} \otimes (dr - r_j^{3-n} \, d\hat{u}) \big ) = 0,\] 
it follows that 
\[\big | \text{\rm tr}_{\partial M^{(j)}}(D^2 \hat{u}) - \text{\rm tr}_{\partial M^{(j)}}(D_\gamma^2 u) \big | \leq o(r_j^{-2})\] 
at each point on $\partial M^{(j)}$. Consequently, 
\begin{equation} 
\label{second.term} 
\big | \text{\rm tr}_{\partial M^{(j)}}(D^2 \hat{u}) - r^{-2} \, \Delta_\gamma u \big | \leq o(r_j^{-2}) 
\end{equation}
at each point on $\partial M^{(j)}$. Combining (\ref{first.term}) and (\ref{second.term}) gives 
\begin{align*} 
&\Big | \text{\rm tr}_{\partial M^{(j)}}(D^2 r) - r_j^{3-n} \, \text{\rm tr}_{\partial M^{(j)}}(D^2 \hat{u}) \\ 
&- (n-1) \, r + \frac{n}{2} \, r^{1-n} \, \text{\rm tr}_\gamma(Q) + r_j^{3-n} \, r^{-2} \, \Delta_\gamma u \Big | \leq o(r_j^{1-n}) 
\end{align*}
at each point on $\partial M^{(j)}$. Using (\ref{pde.for.u}), we conclude that 
\begin{equation}
\label{boundary.trace.of.Hessian}
\Big | \text{\rm tr}_{\partial M^{(j)}}(D^2 r) - r_j^{3-n} \, \text{\rm tr}_{\partial M^{(j)}}(D^2 \hat{u}) - (n-1) \, r - r^{1-n} \, \mu \Big | \leq o(r_j^{1-n}) 
\end{equation}
at each point on $\partial M^{(j)}$. 

We next observe that 
\[|dr|_{\hat{g}}^2 = r^2, \quad \langle dr,d\hat{u} \rangle_{\hat{g}} = 0, \quad |d\hat{u}|_{\hat{g}}^2 \leq O(r^{-2}).\] 
Since $|g-\hat{g}|_{\hat{g}} \leq o(r^{-n})$, it follows that 
\[|dr|_g^2 = r^2 + o(r^{2-n}), \quad \langle dr,d\hat{u} \rangle_g = o(r^{-n}), \quad |d\hat{u}|_g^2 \leq O(r^{-2}).\] 
This implies 
\[\big | |dr - r_j^{3-n} \, d\hat{u}|_g^2 - r^2 \big | \leq o(r_j^{2-n})\] 
at each point on $\partial M^{(j)}$. Consequently, 
\begin{equation} 
\label{norm.of.gradient} 
\big | |dr - r_j^{3-n} \, d\hat{u}|_g - r \big | \leq o(r_j^{1-n}) 
\end{equation} 
at each point on $\partial M^{(j)}$. Using (\ref{boundary.trace.of.Hessian}) and (\ref{norm.of.gradient}), we obtain 
\begin{equation} 
\label{quotient}
\Big | \frac{\text{\rm tr}_{\partial M^{(j)}}(D^2 r) - r_j^{3-n} \, \text{\rm tr}_{\partial M^{(j)}}(D^2 \hat{u})}{|dr - r_j^{3-n} \, d\hat{u}|_g} - (n-1) - r^{-n} \, \mu \Big | \leq o(r_j^{1-n}) 
\end{equation}
at each point on $\partial M^{(j)}$. The assertion follows by combining (\ref{mean.curvature}) and (\ref{quotient}). This completes the proof of Proposition \ref{asymptotic.formula.for.mean.curvature}. \\

After these preparations, we now describe the proof of Theorem \ref{positive.energy.theorem}. Let $\Theta_0$ denote the pull-back of the volume form on $S^1$ under the canonical projection from $S^1 \times T^{n-2}$ to $S^1$. Note that $\Theta_0$ is a closed one-form on $S^1 \times T^{n-2}$. Let $\sigma$ denote the length of the shortest closed curve $\alpha$ in $(S^1 \times T^{n-2},\gamma)$ satisfying $\int_\alpha \Theta_0 \neq 0$. 

We extend $\Theta_0$ to a closed one-form on $[r_0,\infty) \times S^1 \times T^{n-2}$. We may view $\Theta_0$ as a closed one-form which is defined on $M \setminus E$. Let $\theta_0: [r_0,\infty) \times S^1 \times T^{n-2} \to S^1$ denote the canonical projection to the second factor. We may view $\theta_0$ as a map from $M \setminus E$ to $S^1$. The pull-back of the volume form on $S^1$ under the map $\theta_0: M \setminus E \to S^1$ is given by $\Theta_0$.

Let $(\theta_1,\hdots,\theta_{n-2}): [r_0,\infty) \times S^1 \times T^{n-2} \to T^{n-2}$ denote the canonical projection to the third factor. We may view $(\theta_1,\hdots,\theta_{n-2})$ as a map from $M \setminus E$ to $T^{n-2}$. By assumption, the map $(\theta_1,\hdots,\theta_{n-2})$ extends to a smooth map from $M$ to $T^{n-2}$. 

To prove Theorem \ref{positive.energy.theorem}, we argue by contradiction. Suppose that 
\begin{equation} 
\label{mass}
\int_{S^1 \times T^{n-2}} \Big ( n \, \text{\rm tr}_\gamma(Q) + \Big ( \frac{4\pi}{n\sigma} \Big )^n \Big ) \, d\text{\rm vol}_\gamma < 0. 
\end{equation}
Combining (\ref{definition.of.mu}) and (\ref{mass}), we obtain 
\begin{equation} 
\label{lower.bound.for.mu}
2 \, \sigma^n \, \mu > \Big ( \frac{4\pi}{n} \Big )^n. 
\end{equation}
In particular, $\mu > 0$. 

In the following, we assume that $j$ is sufficiently large. We consider the domain $M^{(j)} = M \setminus \{r > r_j + r_j^{3-n} \, \hat{u}\}$. With a suitable choice of orientation, the map $(\theta_0|_{\partial M^{(j)}},\theta_1|_{\partial M^{(j)}},\hdots,\theta_{n-2}|_{\partial M^{(j)}}): \partial M^{(j)} \to S^1 \times T^{n-2}$ has degree $1$. Using Proposition \ref{asymptotic.formula.for.mean.curvature}, we obtain 
\begin{equation} 
\label{lower.bound.for.mean.curvature}
\liminf_{j \to \infty} \Big ( r_j^n \, \inf_{\partial M^{(j)}} (H_{\partial M^{(j)}} - (n-1)) \Big ) \geq \mu. 
\end{equation}
Let $\sigma_j$ denote the length of the shortest closed curve $\alpha$ in $(\partial M^{(j)},g)$ satisfying $\int_\alpha \Theta_0 \neq 0$. Then 
\begin{equation} 
\label{lower.bound.for.systole}
\liminf_{j \to \infty} r_j^{-1} \, \sigma_j \geq \sigma. 
\end{equation} 
Combining (\ref{lower.bound.for.mean.curvature}) and (\ref{lower.bound.for.systole}) gives 
\begin{equation} 
\label{liminf}
\liminf_{j \to \infty} \Big ( 2 \, \sigma_j^n \, \inf_{\partial M^{(j)}} (H_{\partial M^{(j)}} - (n-1)) \Big ) \geq 2 \, \sigma^n \, \mu. 
\end{equation} 
On the other hand, Theorem \ref{higher.dim.inequality.2} implies that 
\begin{equation} 
\label{systolic.inequality}
2 \, \sigma_j^n \, \inf_{\partial M^{(j)}} (H_{\partial M^{(j)}} - (n-1)) \leq \Big ( \frac{4\pi}{n} \Big )^n 
\end{equation}
for each $j$. Combining (\ref{liminf}) and (\ref{systolic.inequality}), we conclude that 
\[2 \, \sigma^n \, \mu \leq \Big ( \frac{4\pi}{n} \Big )^n.\] 
This contradicts (\ref{lower.bound.for.mu}). This completes the proof of Theorem \ref{positive.energy.theorem}. 

\appendix

\section{The second variation formula for weighted area}

In this section, we derive the stability inequality for free boundary minimal hypersurfaces with respect to a conformally modified metric. 

\begin{theorem}
\label{second.variation}
Let $M$ be a compact, orientable manifold of dimension $n$ with non-empty boundary $\partial M$. Let $g$ be a Riemannian metric on $M$, and let $\rho$ be a smooth positive function on $M$. Suppose that $\Sigma$ is an orientable hypersurface in $M$ such that $\partial \Sigma \subset \partial M$ and $\Sigma$ meets $\partial M$ orthogonally along $\partial \Sigma$. If $\Sigma$ is a stable free boundary minimal hypersurface in $(M,\rho^{\frac{2}{n-1}} \, g)$, then 
\begin{align*} 
&\int_\Sigma \rho \, |\nabla^\Sigma \zeta|^2 - \int_\Sigma \rho \, \text{\rm Ric}_M(\nu_\Sigma,\nu_\Sigma) \, \zeta^2 - \int_\Sigma \rho \, |h_\Sigma|^2 \, \zeta^2 \\ 
&+ \int_\Sigma (D_M^2 \rho)(\nu_\Sigma,\nu_\Sigma) \, \zeta^2 - \int_\Sigma \rho^{-1} \, \langle \nabla^M \rho,\nu_\Sigma \rangle^2 \, \zeta^2 \\ 
&- \int_{\partial \Sigma} \rho \, h_{\partial M}(\nu_\Sigma,\nu_\Sigma) \, \zeta^2 \geq 0 
\end{align*}
for every test function $\zeta \in C^\infty(\Sigma)$. 
\end{theorem}

\textbf{Proof.} 
Let $\eta$ denote the outward-pointing unit normal vector field to $\partial M$ with respect to $g$, and let $h_{\partial M}$ denote the second fundamental form of $\partial M$ with respect to $g$. We consider the conformal metric $\tilde{g} = \rho^{\frac{2}{n-1}} \, g$. The unit normal vector field to $\partial M$ with respect to $\tilde{g}$ is given by $\tilde{\eta} = \rho^{-\frac{1}{n-1}} \, \eta$, and the second fundamental form of $\partial M$ with respect to $\tilde{g}$ is given by 
\[\tilde{h}_{\partial M} = \rho^{\frac{1}{n-1}} \, \Big ( h_{\partial M} + \frac{1}{n-1} \, \rho^{-1} \, d\rho(\eta) \, g \Big ).\] 
The unit normal vector field to $\Sigma$ with respect to $\tilde{g}$ is given by $\tilde{\nu}_\Sigma = \rho^{-\frac{1}{n-1}} \, \nu_\Sigma$, and the second fundamental form of $\Sigma$ with respect to $\tilde{g}$ is given by 
\[\tilde{h}_\Sigma = \rho^{\frac{1}{n-1}} \, \Big ( h_\Sigma + \frac{1}{n-1} \, \rho^{-1} \, d\rho(\nu_\Sigma) \, g \Big ).\] 
Since $\Sigma$ is a minimal hypersurface in $(M,\tilde{g})$, it follows that $H_\Sigma + \rho^{-1} \, d\rho(\nu_\Sigma) = 0$. This implies 
\[|\tilde{h}_\Sigma|_{\tilde{g}}^2 = \rho^{-\frac{2}{n-1}} \, \Big ( |h_\Sigma|^2 - \frac{1}{n-1} \, \rho^{-2} \, (d\rho(\nu_\Sigma))^2 \Big ).\] 
The Ricci tensor of $\tilde{g}$ is related to the Ricci tensor of $g$ by the formula 
\begin{align*} 
\text{\rm Ric}_{\tilde{g}} 
&= \text{\rm Ric} - \Big ( 1- \frac{1}{n-1} \Big ) \, \rho^{-1} \, D^2 \rho - \frac{1}{n-1} \, \rho^{-1} \, \Delta \rho \, g \\ 
&+ \Big ( 1 - \frac{1}{(n-1)^2} \Big ) \, \rho^{-2} \, d\rho \otimes d\rho + \frac{1}{(n-1)^2} \, \rho^{-2} \, |d\rho|^2 \, g 
\end{align*} 
(see \cite{Besse}, Theorem 1.159). Since $\Sigma$ is a stable free boundary minimal hypersurface in $(M,\tilde{g})$, the second variation formula implies that 
\begin{align*} 
&\int_\Sigma |d\zeta|_{\tilde{g}}^2 \, d\text{\rm vol}_{\tilde{g}} - \int_\Sigma \text{\rm Ric}_{\tilde{g}}(\tilde{\nu}_\Sigma,\tilde{\nu}_\Sigma) \, \zeta^2 \, d\text{\rm vol}_{\tilde{g}} - \int_\Sigma |\tilde{h}_\Sigma|_{\tilde{g}}^2 \, \zeta^2 \, d\text{\rm vol}_{\tilde{g}} \\ 
&- \int_{\partial \Sigma} \tilde{h}_{\partial M}(\tilde{\nu}_\Sigma,\tilde{\nu}_\Sigma) \, \zeta^2 \, d\text{\rm vol}_{\tilde{g}} \geq 0 
\end{align*}
for every test function $\zeta \in C^\infty(\Sigma)$ (see \cite{Ambrozio-Carlotto-Sharp}). In the next step, we replace $\zeta$ by $\rho^{\frac{1}{n-1}} \, \zeta$. Moreover, we convert all the geometric quantities back to the metric $g$. This gives 
\begin{align*} 
&\int_\Sigma \rho^{1-\frac{2}{n-1}} \, |\nabla^\Sigma (\rho^{\frac{1}{n-1}} \, \zeta)|^2 - \int_\Sigma \rho \, \text{\rm Ric}_M(\nu_\Sigma,\nu_\Sigma) \, \zeta^2 - \int_\Sigma \rho \, |h_\Sigma|^2 \, \zeta^2 \\ 
&+ \Big ( 1 - \frac{1}{n-1} \Big ) \int_\Sigma (D_M^2 \rho)(\nu_\Sigma,\nu_\Sigma) \, \zeta^2 + \frac{1}{n-1} \int_\Sigma \Delta_M \rho \, \zeta^2 \\ 
&- \Big ( 1 - \frac{1}{n-1} - \frac{1}{(n-1)^2} \Big ) \int_\Sigma \rho^{-1} \, \langle \nabla^M \rho,\nu_\Sigma \rangle^2 \, \zeta^2 \\ 
&- \frac{1}{(n-1)^2} \int_\Sigma \rho^{-1} \, |\nabla^M \rho|^2 \, \zeta^2 \\ 
&- \int_{\partial \Sigma} \rho \, h_{\partial M}(\nu_\Sigma,\nu_\Sigma) \, \zeta^2 - \frac{1}{n-1} \int_{\partial \Sigma} \langle \nabla^\Sigma \rho,\eta \rangle \, \zeta^2 \geq 0 
\end{align*}
for every test function $\zeta \in C^\infty(\Sigma)$. Using the divergence theorem, we obtain 
\begin{align*} 
&\int_\Sigma \rho^{1-\frac{2}{n-1}} \, |\nabla^\Sigma (\rho^{\frac{1}{n-1}} \, \zeta)|^2 - \frac{1}{n-1} \int_{\partial \Sigma} \langle \nabla^\Sigma \rho,\eta \rangle \, \zeta^2 \\ 
&= \int_\Sigma \rho^{1-\frac{2}{n-1}} \, |\nabla^\Sigma (\rho^{\frac{1}{n-1}} \, \zeta)|^2 - \frac{1}{n-1} \int_\Sigma \text{\rm div}_\Sigma(\zeta^2 \, \nabla^\Sigma \rho) \\ 
&= \int_\Sigma \rho \, |\nabla^\Sigma \zeta|^2 - \frac{1}{n-1} \int_\Sigma \Delta_\Sigma \rho \, \zeta^2 + \frac{1}{(n-1)^2} \int_\Sigma \rho^{-1} \, |\nabla^\Sigma \rho|^2 \, \zeta^2 
\end{align*}
for every test function $\zeta \in C^\infty(\Sigma)$. Putting these facts together, we conclude that 
\begin{align*} 
&\int_\Sigma \rho \, |\nabla^\Sigma \zeta|^2 - \int_\Sigma \rho \, \text{\rm Ric}_M(\nu_\Sigma,\nu_\Sigma) \, \zeta^2 - \int_\Sigma \rho \, |h_\Sigma|^2 \, \zeta^2 \\ 
&+ \int_\Sigma (D_M^2 \rho)(\nu_\Sigma,\nu_\Sigma) \, \zeta^2 + \frac{1}{n-1} \int_\Sigma (\Delta_M \rho - (D_M^2 \rho)(\nu_\Sigma,\nu_\Sigma) - \Delta_\Sigma \rho) \, \zeta^2 \\ 
&- \Big ( 1 - \frac{1}{n-1} - \frac{1}{(n-1)^2} \Big ) \int_\Sigma \rho^{-1} \, \langle \nabla^M \rho,\nu_\Sigma \rangle^2 \, \zeta^2 \\ 
&- \frac{1}{(n-1)^2} \int_\Sigma \rho^{-1} \, (|\nabla^M \rho|^2 - |\nabla^\Sigma \rho|^2) \, \zeta^2 \\ 
&- \int_{\partial \Sigma} \rho \, h_{\partial M}(\nu_\Sigma,\nu_\Sigma) \, \zeta^2 \geq 0 
\end{align*}
for every test function $\zeta \in C^\infty(\Sigma)$. The assertion now follows from the identities 
\[\Delta_M \rho - (D_M^2 \rho)(\nu_\Sigma,\nu_\Sigma) - \Delta_\Sigma \rho = H_\Sigma \, \langle \nabla^M \rho,\nu_\Sigma \rangle = -\rho^{-1} \, \langle \nabla^M \rho,\nu_\Sigma \rangle^2\] 
and 
\[|\nabla^M \rho|^2 - |\nabla^\Sigma \rho|^2 = \langle \nabla^M \rho,\nu_\Sigma \rangle^2.\]
This completes the proof of Theorem \ref{second.variation}. 

\section{Existence and regularity of free boundary minimal hypersurfaces}

In this section, we recall some well known results concerning the existence and regularity of free boundary minimal hypersurfaces.

\begin{theorem}[cf. H.~Federer \cite{Federer}; M.~Gr\"uter \cite{Grueter}]
\label{existence.and.regularity} 
Let us fix an integer $3 \leq n \leq 7$. Let $M$ be a compact, orientable manifold of dimension $n$ with non-empty boundary $\partial M$. Let $g$ be a Riemannian metric on $M$. We assume that the mean curvature of $\partial M$ with respect to $g$ is strictly positive. Let $\Omega$ be a closed $(n-2)$-form on $\partial M$. Let $\tilde{\Sigma}$ be a compact, embedded, orientable hypersurface in $M$ such that $\partial \tilde{\Sigma} \subset \partial M$ and $\int_{\partial \tilde{\Sigma}} \Omega \neq 0$. Then we can find a compact, connected, embedded, orientable hypersurface $\Sigma$ with the following properties: 
\begin{itemize} 
\item The boundary $\partial \Sigma$ is contained in $\partial M$. Moreover, $\Sigma$ meets $\partial M$ orthogonally along $\partial \Sigma$. 
\item $\Sigma$ is a stable free boundary minimal hypersurface.
\item $\int_{\partial \Sigma} \Omega \neq 0$. 
\end{itemize}
\end{theorem}

In the remainder of this section, we explain how Theorem \ref{existence.and.regularity} follows from results of Federer \cite{Federer} and Gr\"uter \cite{Grueter}. Let $\hat{M}$ be a compact manifold which contains the given manifold $M$ in its interior. We extend the given Riemannian metric $g$ on $M$ to a Riemannian metric on $\hat{M}$. If $\varepsilon_0 > 0$ is sufficiently small, then the set $\{x \in \hat{M}: d_{\hat{M}}(x,\partial M) \leq \varepsilon_0\}$ can be identified with $\partial M \times [-\varepsilon_0,\varepsilon_0]$ via the normal exponential map. For $\varepsilon \in (0,\varepsilon_0]$, we define $\hat{M}_\varepsilon = \{x \in \hat{M}: d_{\hat{M}}(x,M) \leq \varepsilon\}$. If we choose $\varepsilon_0 > 0$ sufficiently small, then the boundary $\partial \hat{M}_\varepsilon$ is strictly mean convex for each $\varepsilon \in (0,\varepsilon_0]$. We denote by $f: \hat{M}_{\varepsilon_0} \to M$ the nearest point projection from $\hat{M}_{\varepsilon_0}$ to $M$. Note that $f$ is Lipschitz continuous and $f|_M = \text{\rm id}$. For each $\varepsilon \in (0,\varepsilon_0]$, the restriction $f|_{\partial \hat{M}_\varepsilon}$ is a volume-decreasing map from $\partial \hat{M}_\varepsilon$ to $\partial M$. 

\begin{lemma} 
\label{differential.of.f}
Let $x \in \hat{M}_{\varepsilon_0} \setminus M$. Then 
\[|df_x(v_1) \wedge \hdots \wedge df_x(v_{n-1})| \leq |v_1 \wedge \hdots \wedge v_{n-1}|\] 
for all tangent vectors $v_1,\hdots,v_{n-1} \in T_x \hat{M}$.
\end{lemma} 

\textbf{Proof.} 
Let us fix a real number $\varepsilon \in (0,\varepsilon_0]$, a point $x \in \partial \hat{M}_\varepsilon$, and tangent vectors $v_1,\hdots,v_{n-1} \in T_x \hat{M}$. For each $1 \leq i \leq n-1$, we denote by $w_i$ the orthogonal projection of $v_i$ to the tangent space $T_x(\partial \hat{M}_\varepsilon)$. Since the restriction $f|_{\partial \hat{M}_\varepsilon}$ is a volume-decreasing map from $\partial \hat{M}_\varepsilon$ to $\partial M$, we know that 
\begin{equation} 
\label{step.1}
|df_x(w_1) \wedge \hdots \wedge df_x(w_{n-1})| \leq |w_1 \wedge \hdots \wedge w_{n-1}|. 
\end{equation}
Since the matrix $\{\langle v_i,v_j \rangle - \langle w_i,w_j \rangle\}_{1 \leq i,j \leq n-1}$ is weakly positive definite, we obtain 
\[\det \{\langle w_i,w_j \rangle\}_{1 \leq i,j \leq n-1} \leq \det \{\langle v_i,v_j \rangle\}_{1 \leq i,j \leq n-1},\] 
hence 
\begin{equation} 
\label{step.2}
|w_1 \wedge \hdots \wedge w_{n-1}| \leq |v_1 \wedge \hdots \wedge v_{n-1}|. 
\end{equation}
Finally, it follows from the definition of $w_i$ that $df_x(w_i) = df_x(v_i)$ for $1 \leq i \leq n-1$. This implies 
\begin{equation} 
\label{step.3} 
|df_x(w_1) \wedge \hdots \wedge df_x(w_{n-1})| = |df_x(v_1) \wedge \hdots \wedge df_x(v_{n-1})|. 
\end{equation} 
Combining (\ref{step.1}), (\ref{step.2}), and (\ref{step.3}), the assertion follows. This completes the proof of Lemma \ref{differential.of.f}. \\

By the Nash embedding theorem, the manifold $(\hat{M},g)$ can be isometrically embedded into $\mathbb{R}^N$ for some large integer $N$. Let $\delta$ be a positive real number, and let $U \subset \mathbb{R}^N$ denote the set of all points in $\mathbb{R}^N$ that have distance less than $\delta$ from $M$. If $\delta>0$ is chosen sufficiently small (depending on $\varepsilon_0$), then the nearest point projection gives a smooth map $\pi: U \to \hat{M}_{\varepsilon_0}$. We define a map $F: U \to M$ by $F = f \circ \pi$. Note that $F$ is Lipschitz continuous and $F|_M = \text{\rm id}$. 

We next recall the definition of an integer multiplicity rectifiable current from Leon Simon's notes \cite{Simon}. 

\begin{definition}[cf. L.~Simon \cite{Simon}, Chapter 6, \S 3]
\label{integer.multiplicity.rectifiable.currents}
Let $m$ be an integer with $1 \leq m \leq N$, and let $T$ be an $m$-dimensional current in $\mathbb{R}^N$. We say that $T$ is an integer multiplicity rectifiable current if we can find a set $A \subset \mathbb{R}^N$, a positive function $\theta: A \to \mathbb{Z}$, and a function $\xi: A \to \Lambda^m(\mathbb{R}^N)$ with the following properties:
\begin{itemize}
\item The set $A$ is $\mathcal{H}^m$-measurable and countably $m$-rectifiable. 
\item The function $\theta: A \to \mathbb{Z}$ is locally $\mathcal{H}^m$-integrable. 
\item The function $\xi: A \to \Lambda^m(\mathbb{R}^N)$ is $\mathcal{H}^m$-measurable. 
\item For $\mathcal{H}^m$-a.e. point $x \in A$, $\xi(x)$ can be expressed in the form $v_1 \wedge \hdots \wedge v_m$, where $v_1,\hdots,v_m$ form an orthonormal basis for the approximate tangent space $T_x A$.
\item We have 
\[T(\omega) = \int_A \langle \omega(x),\xi(x) \rangle \, \theta(x) \, d\mathcal{H}^m(x)\] 
for every smooth $m$-form $\omega$ on $\mathbb{R}^N$ with compact support.
\end{itemize}
\end{definition}

If $T$ is an $m$-dimensional integer multiplicity rectifiable current with compact support, then the mass of $T$ is given by 
\[\text{\rm \bf M}(T) = \int_A \theta(x) \, d\mathcal{H}^m(x)\] 
(see \cite{Simon}, Chapter 6, \S 3). 

\begin{lemma}
\label{area.decreasing.property}
Suppose that $T$ is an $(n-1)$-dimensional integer multiplicity rectifiable current satisfying $\text{\rm supp}(T) \subset \hat{M}_{\varepsilon_0} \cap U$. Then $F_\#(T)$ is an $(n-1)$-dimensional integer multiplicity rectifiable current and $\text{\rm \bf M}(F_\#(T)) \leq \text{\rm \bf M}(T)$.
\end{lemma} 

\textbf{Proof.} 
We use the notation from Definition \ref{integer.multiplicity.rectifiable.currents}. Let $J_F^A$ denote the unoriented Jacobian as defined in Leon Simon's notes (see \cite{Simon}, Chapter 3, \S 2). Note that $J_F^A(x)$ is defined for $\mathcal{H}^{n-1}$-a.e. point $x \in A$, and $J_F^A(x)$ is nonnegative. Since $\text{\rm supp}(T) \subset \hat{M}_{\varepsilon_0} \cap U$, the approximate tangent space $T_x A$ is contained in $T_x \hat{M}$ for $\mathcal{H}^{n-1}$-a.e. point $x \in A$. Since $F|_M = \text{\rm id}$, we obtain $J_F^A(x) \leq 1$ for $\mathcal{H}^{n-1}$-a.e. point $x \in A \cap M$. Using Lemma \ref{differential.of.f}, we obtain $J_F^A(x) \leq 1$ for $\mathcal{H}^{n-1}$-a.e. point $x \in A \setminus M$. Thus, $J_F^A(x) \leq 1$ for $\mathcal{H}^{n-1}$-a.e. point $x \in A$. 

It is shown in Leon Simon's notes that $F_\#(T)$ is an $(n-1)$-dimensional integer multiplicity rectifiable current (see \cite{Simon}, Chapter 6, \S 3). The mass of $F_\#(T)$ satisfies the estimate 
\[\text{\rm \bf M}(F_\#(T)) \leq \int_{F(A)} \mathcal{N}(y) \, d\mathcal{H}^{n-1}(y)\] 
(see \cite{Simon}, Chapter 6, \S 3), where the function $\mathcal{N}: F(A) \to \mathbb{R}$ is defined by 
\[\mathcal{N}(y) = \sum_{x \in A, \, F(x)=y, \, J_F^A(x) > 0} \theta(x)\] 
for $y \in F(A)$. This implies 
\[\text{\rm \bf M}(F_\#(T)) \leq \int_A J_F^A(x) \, \theta(x) \, d\mathcal{H}^{n-1}(x).\] 
Since $J_F^A(x) \leq 1$ for $\mathcal{H}^{n-1}$-a.e. point $x \in A$, we conclude that 
\[\text{\rm \bf M}(F_\#(T)) \leq \int_A \theta(x) \, d\mathcal{H}^{n-1}(x) = \text{\rm \bf M}(T).\] 
This completes the proof of Lemma \ref{area.decreasing.property}. \\

In the following, we recall some basic definitions from Federer's work \cite{Federer}. Given a positive integer $m$ and a compact set $K \subset \mathbb{R}^N$, let $\mathcal{R}_{m,K}(\mathbb{R}^N)$ be defined as in \cite{Federer}, Section 4.1.24. It follows directly from the definition of $\mathcal{R}_{m,K}(\mathbb{R}^N)$ that $\text{\rm supp}(T) \subset K$ for all $T \in \mathcal{R}_{m,K}(\mathbb{R}^N)$. For every positive integer $m$, we define 
\[\mathcal{R}_m(\mathbb{R}^N) = \bigcup_K \mathcal{R}_{m,K}(\mathbb{R}^N),\] 
where the union is taken over all compact sets $K \subset \mathbb{R}^N$ (see \cite{Federer}, Section 4.1.24). In particular, if $T \in \mathcal{R}_m(\mathbb{R}^N)$, then $T$ has compact support.

\begin{proposition}[cf. H.~Federer \cite{Federer}]
\label{characterization.of.R}
Let $T$ be an $m$-dimensional current in $\mathbb{R}^N$ with compact support. Then $T \in \mathcal{R}_m(\mathbb{R}^N)$ if and only if $T$ is an integer multiplicity rectifiable current in the sense of Definition \ref{integer.multiplicity.rectifiable.currents}.
\end{proposition} 

\textbf{Proof.} 
This follows from Theorem 4.1.28 in \cite{Federer}. We specifically use the equivalence of statements (1) and (4). This completes the proof of Proposition \ref{characterization.of.R}. \\

Let $m$ be an integer with $1 \leq m \leq N$. The space of $m$-dimensional integral flat chains on $\mathbb{R}^N$ is defined by  
\[\mathcal{F}_m(\mathbb{R}^N) = \{P+\partial Q: P \in \mathcal{R}_m(\mathbb{R}^N), \, Q \in \mathcal{R}_{m+1}(\mathbb{R}^N)\}\] 
(see \cite{Federer}, Section 4.1.24). Note that $\mathcal{R}_m(\mathbb{R}^N) \subset \mathcal{F}_m(\mathbb{R}^N)$. The space of $m$-dimensional integral flat cycles is defined by 
\[\mathcal{Z}_m(M,\partial M) = \{T \in \mathcal{F}_m(\mathbb{R}^N): \text{\rm supp}(T) \subset M, \, \text{\rm supp}(\partial T) \subset \partial M\}\] 
(see \cite{Federer}, Section 4.4.1). The space of $m$-dimensional integral flat boundaries is defined by 
\begin{align*} 
&\mathcal{B}_m(M,\partial M) \\ 
&= \{P + \partial Q: P \in \mathcal{F}_m(\mathbb{R}^N), \, Q \in \mathcal{F}_{m+1}(\mathbb{R}^N), \, \text{\rm supp}(P) \subset \partial M, \, \text{\rm supp}(Q) \subset M\} 
\end{align*} 
(see \cite{Federer}, Section 4.4.1). Note that $\mathcal{B}_m(M,\partial M) \subset \mathcal{Z}_m(M,\partial M)$. \\

We now continue the proof of Theorem \ref{existence.and.regularity}. By assumption, $\Omega$ is a closed $(n-2)$-form on $\partial M$. We may extend $\Omega$ to an $(n-2)$-form defined on $\mathbb{R}^N$ such that $\Omega$ has compact support and $d\Omega$ vanishes in an open neighborhood of $\partial M$. (To construct this extension, we map a small tubular neighborhood of $\partial M$ in $\mathbb{R}^N$ to $\partial M$ using the nearest point projection. We then consider the pull-back of $\Omega$ under this map.) 

\begin{lemma}
\label{integral.of.T.over.Omega}
If $T \in \mathcal{B}_{n-1}(M,\partial M)$, then $\partial T(\Omega) = 0$.
\end{lemma}

\textbf{Proof.} 
Let $T \in \mathcal{B}_{n-1}(M,\partial M)$. We may write $T = P + \partial Q$, where $P \in \mathcal{F}_m(\mathbb{R}^N)$, $Q \in \mathcal{F}_{m+1}(\mathbb{R}^N)$, $\text{\rm supp}(P) \subset \partial M$, and $\text{\rm supp}(Q) \subset M$. This implies 
\[\partial T(\Omega) = \partial P(\Omega) = P(d\Omega) = 0.\]
In the last step, we have used the fact that $\text{\rm supp}(P) \subset \partial M$ and $d\Omega$ vanishes in an open neighborhood of $\partial M$. This completes the proof of Lemma \ref{integral.of.T.over.Omega}. \\

Let $\tilde{\Sigma}$ denote the hypersurface in Theorem \ref{existence.and.regularity}. The hypersurface $\tilde{\Sigma}$ defines a current $\tilde{S} \in \mathcal{Z}_{n-1}(M,\partial M)$ such that $\partial \tilde{S}(\Omega) \neq 0$. It follows from results in Section 5.1.6 in \cite{Federer} that we can find a current $S \in \mathcal{R}_{n-1}(\mathbb{R}^N)$ with the following properties: 
\begin{itemize}
\item We have $S-\tilde{S} \in \mathcal{B}_{n-1}(M,\partial M)$. In particular, $S \in \mathcal{Z}_{n-1}(M,\partial M)$.
\item The current $S$ is homologically area-minimizing with respect to $(M,\partial M)$. In other words, $\text{\rm \bf M}(S) \leq \text{\rm \bf M}(S+X)$ for all currents $X \in \mathcal{B}_{n-1}(M,\partial M) \cap \mathcal{R}_{n-1}(\mathbb{R}^N)$.
\end{itemize} 
Since $S \in \mathcal{R}_{n-1}(\mathbb{R}^N)$, $S$ is an $(n-1)$-dimensional integer multiplicity rectifiable current by Proposition \ref{characterization.of.R}. Since $S-\tilde{S} \in \mathcal{B}_{n-1}(M,\partial M)$, Lemma \ref{integral.of.T.over.Omega} implies that $\partial S(\Omega) = \partial \tilde{S}(\Omega)$. In particular, $\partial S(\Omega) \neq 0$.

\begin{lemma}
\label{S.locally.area.minimizing}
Let $p$ be an arbitrary point in $M$. Then we can find an open set $O \subset \mathbb{R}^N$ containing $p$ such that $\text{\rm \bf M}(S) \leq \text{\rm \bf M}(S+X)$ for all $(n-1)$-dimensional integer multiplicity rectifiable currents $X$ satisfying $\text{\rm supp}(X) \subset M \cap O$ and $\text{\rm supp}(\partial X) \subset \partial M \cap O$. 
\end{lemma} 

\textbf{Proof.}
Since $S$ is homologically area-minimizing with respect to $(M,\partial M)$, the results in Section 5.1.6 in \cite{Federer} imply that $S$ is locally area-minimizing with respect to $(M,\partial M)$. Consequently, we can find an open set $O \subset \mathbb{R}^N$ containing $p$ such that $\text{\rm \bf M}(S) \leq \text{\rm \bf M}(S+X)$ for all currents $X \in \mathcal{R}_{n-1}(\mathbb{R}^N)$ satisfying $\text{\rm supp}(X) \subset M \cap O$ and $\text{\rm supp}(\partial X) \subset \partial M \cap O$. Finally, every $(n-1)$-dimensional integer multiplicity rectifiable current with compact support belongs to $\mathcal{R}_{n-1}(\mathbb{R}^N)$ by Proposition \ref{characterization.of.R}. This completes the proof of Lemma \ref{S.locally.area.minimizing}. \\

\begin{lemma}
\label{S.locally.area.minimizing.stronger.version}
Let $p$ be an arbitrary point in $M$. Then we can find an open set $\hat{O} \subset \mathbb{R}^N$ containing $p$ such that $\text{\rm \bf M}(S) \leq \text{\rm \bf M}(S+X)$ for all $(n-1)$-dimensional integer multiplicity rectifiable currents $X$ satisfying $\text{\rm supp}(X) \subset \hat{M}_{\varepsilon_0} \cap \hat{O}$ and $\text{\rm supp}(\partial X) \subset \partial M \cap \hat{O}$. 
\end{lemma}

\textbf{Proof.} 
By Lemma \ref{S.locally.area.minimizing}, we can find an open set $O \subset \mathbb{R}^N$ containing $p$ such that 
\[\text{\rm \bf M}(S) \leq \text{\rm \bf M}(S+X)\] 
for all $(n-1)$-dimensional integer multiplicity rectifiable currents $X$ satisfying $\text{\rm supp}(X) \subset M \cap O$ and $\text{\rm supp}(\partial X) \subset \partial M \cap O$. We define 
\[\hat{O} = \{x \in U: F(x) \in O\}.\] 
Clearly, $\hat{O} \subset \mathbb{R}^N$ is an open set containing $p$.

We claim that $\hat{O}$ has the desired property. To see this, suppose that $X$ is an $(n-1)$-dimensional integer multiplicity rectifiable current satisfying $\text{\rm supp}(X) \subset \hat{M}_{\varepsilon_0} \cap \hat{O}$ and $\text{\rm supp}(\partial X) \subset \partial M \cap \hat{O}$. Since $\hat{O} \subset U$, it follows from Proposition \ref{area.decreasing.property} that $F_\#(X)$ is an $(n-1)$-dimensional integer multiplicity rectifiable current. The results in Section 4.1.14 in \cite{Federer} imply that 
\[\text{\rm supp}(F_\#(X)) \subset F(\text{\rm supp}(X)) \subset F(\hat{O})\] 
and 
\[\text{\rm supp}(\partial (F_\#(X))) = \text{\rm supp}(F_\#(\partial X)) \subset F(\text{\rm supp}(\partial X)) \subset F(\partial M \cap \hat{O}).\] 
Note that $F(\hat{O}) \subset M \cap O$ by definition of $\hat{O}$. Moreover, since $F|_M = \text{\rm id}$, we obtain $F(\partial M \cap \hat{O}) \subset \partial M \cap O$. Therefore, 
\[\text{\rm supp}(F_\#(X)) \subset M \cap O\] 
and 
\[\text{\rm supp}(\partial (F_\#(X))) \subset \partial M \cap O.\] 
In view of our choice of $O$, it follows that 
\begin{equation} 
\label{consequence.of.area.minimizing.property}
\text{\rm \bf M}(S) \leq \text{\rm \bf M}(S+F_\#(X)). 
\end{equation}
On the other hand, since $S$ and $X$ are $(n-1)$-dimensional integer multiplicity rectifiable currents, the sum $S+X$ is an $(n-1)$-dimensional integer multiplicity rectifiable current (see \cite{Simon}, Chapter 6, \S 3). Since $\text{\rm supp}(S+X) \subset \hat{M}_{\varepsilon_0} \cap U$, Lemma \ref{area.decreasing.property} implies that 
\begin{equation} 
\label{consequence.of.area.decreasing.property}
\text{\rm \bf M}(F_\#(S+X)) \leq \text{\rm \bf M}(S+X). 
\end{equation}
Finally, since $\text{\rm supp}(S) \subset M$ and $F|_M = \text{\rm id}$, it follows that $F_\#(S) = S$ (see \cite{Federer}, Section 4.1.15). Using (\ref{consequence.of.area.minimizing.property}) and (\ref{consequence.of.area.decreasing.property}), we conclude that 
\[\text{\rm \bf M}(S) \leq \text{\rm \bf M}(S+F_\#(X)) = \text{\rm \bf M}(F_\#(S+X)) \leq \text{\rm \bf M}(S+X).\] 
This completes the proof of Lemma \ref{S.locally.area.minimizing.stronger.version}. \\

Lemma \ref{S.locally.area.minimizing.stronger.version} allows us to apply Gr\"uter's regularity theorem \cite{Grueter} (which was stated in Euclidean space, but extends to the manifold setting). Since $n \leq 7$, Theorem 4.7 in \cite{Grueter} implies that the singular set of $S$ is empty. Since $\partial S(\Omega) \neq 0$, we can find a connected component of $\text{\rm supp}(S)$, denoted by $\Sigma$, with the property that $\int_{\partial \Sigma} \Omega \neq 0$. This completes the proof of Theorem \ref{existence.and.regularity}.

\end{document}